\documentclass[12pt,a4paper]{article}
\usepackage{latexsym,amsfonts,amssymb,amsmath,amscd,graphics,epic,theorem}
\usepackage[all]{xy}
\usepackage{amsmath}
\usepackage{amssymb, amstext}
\baselineskip 18pt \textwidth 15cm  \theoremstyle{plain}
\newtheorem{lemma}{Lemma}[subsection]
\newtheorem{proposition}[lemma]{Proposition}
\newtheorem{remark}[lemma]{Remark}
\newtheorem{example}[lemma]{Example}
\newtheorem{theorem}[lemma]{Theorem}
\newtheorem{definition}{Definition}[subsection]
\newtheorem{notation}[definition]{Notation}
\newtheorem{result}{Result}[subsection]
\newtheorem{property}{}
\newtheorem{corollary}[lemma]{Corollary}

\newtheorem{case}{Case}
{\theorembodyfont{\rmfamily} 
\begin{document}
\newcommand{\pperp}{\hbox{$\perp\hskip-6pt\perp$}}
\newcommand{\N}{{\mathbb N}}
\newcommand{\PP}{{\mathbb P}}
\newcommand{\Z}{{\mathbb Z}}
\newcommand{\Q}{{\mathbb Q}}
\newcommand{\R}{{\mathbb R}}
\newcommand{\C}{{\mathbb C}}
\newcommand{\K}{{\mathbb K}}
\newcommand{\F}{{\mathbb F}}
\newcommand{\proofend}{\hfill$\Box$\smallskip}
\newcommand{\eps}{{\varepsilon}}
\newcommand{\ko}{{\mathcal O}}
\newcommand{\wx}{{\widetilde x}}
\newcommand{\wz}{{\widetilde z}}
\newcommand{\wa}{{\widetilde a}}
\newcommand{\bz}{{\boldsymbol z}}
\newcommand{\bp}{{\boldsymbol p}}
\newcommand{\wy}{{\widetilde y}}
\newcommand{\wc}{{\widetilde c}}
\newcommand{\bi}{{\omega}}
\newcommand{\bx}{{\boldsymbol x}}
\newcommand{\Log}{{\operatorname{Log}}}
\newcommand{\pr}{{\operatorname{pr}}}
\newcommand{\Graph}{{\operatorname{Graph}}}
\newcommand{\jet}{{\operatorname{jet}}}
\newcommand{\Tor}{{\operatorname{Tor}}}
\newcommand{\sqh}{{\operatorname{sqh}}}
\newcommand{\const}{{\operatorname{const}}}
\newcommand{\Arc}{{\operatorname{Arc}}}
\newcommand{\Sing}{{\operatorname{Sing}}}
\newcommand{\Span}{{\operatorname{Span}}}
\newcommand{\Aut}{{\operatorname{Aut}}}
\newcommand{\Ker}{{\operatorname{Ker}}}
\newcommand{\Int}{{\operatorname{Int}}}
\newcommand{\Aff}{{\operatorname{Aff}}}
\newcommand{\Area}{{\operatorname{Area}}}
\newcommand{\val}{{\operatorname{Val}}}
\newcommand{\conv}{{\operatorname{conv}}}
\newcommand{\rk}{{\operatorname{rk}}}
\newcommand{\ow}{{\overline w}}
\newcommand{\ov}{{\overline v}}
\newcommand{\Sc}{{\cal S}}
\newcommand{\G}{{\cal G}}
\newcommand{\T}{{\cal T}}
\newcommand{\red}{{\operatorname{red}}}
\newcommand{\kc}{{\cal C}}
\newcommand{\ki}{{\cal I}}
\newcommand{\kj}{{\cal J}}
\newcommand{\ke}{{\cal E}}
\newcommand{\kz}{{\cal Z}}
\newcommand{\tet}{{\theta}}
\newcommand{\Del}{{\Delta}}
\newcommand{\bet}{{\beta}}
\newcommand{\mm}{{\mathfrak m}}
\newcommand{\kap}{{\kappa}}
\newcommand{\del}{{\delta}}
\newcommand{\sig}{{\sigma}}
\newcommand{\alp}{{\alpha}}
\newcommand{\Sig}{{\Sigma}}
\newcommand{\Gam}{{\Gamma}}
\newcommand{\gam}{{\gamma}}
\newcommand{\Lam}{{\Lambda}}
\newcommand{\lam}{{\lambda}}
\newcommand{\om}{{\omega}}
\newcommand{\cD}{{\mathcal{D}}}
\newcommand{\Fre}{{Fr\'{e}chet \,}}
\title{Schwartz functions on Nash manifolds}
\author{Avraham Aizenbud and Dmitry Gourevitch \thanks{ Avraham Aizenbud and Dmitry Gourevitch, Faculty of Mathematics and Computer
Science, The Weizmann Institute of Science POB 26, Rehovot 76100,
ISRAEL. %Tel: 972-8-934-4209
E-mails: aizenr@yahoo.com, guredim@yahoo.com.\smallskip
 }   }
%\date{}
\maketitle

\begin{abstract} In this paper we extend the notions of Schwartz
functions, tempered functions and generalized Schwartz functions
to Nash (i.e. smooth semi-algebraic) manifolds. % and moreover, to
%the notions of Schwartz sections, tempered sections and
%generalized Schwartz sections of Nash bundles.
We reprove for this
case classically known properties of Schwartz functions on $\R^n$
and build some additional tools %prove some additional properties
which are important in
representation theory.
\end{abstract}

 \tableofcontents
\section{Introduction}

Let us start with the following motivating example. Consider the
circle $S^1$, let $N \subset S^1$ be the north pole and denote
$U:=S^1 \setminus N$. Note that $U$ is diffeomorphic to $\R$ via
the stereographic projection. Consider the space $\cD(S^1)$ of
distributions on $S^1$, that is the space of continuous linear
functionals on the \Fre space $C^{\infty}(S^1)$. Consider the
subspace $\cD_{S^1}(N) \subset \cD(S^1)$ consisting of all
distributions supported at $N$. Then the quotient $ \cD(S^1) /
\cD_{S^1}(N)$ will not be the space of distributions on $U$.
However, it will be the space $\Sc^*(U)$ of Schwartz distributions
on $U$, that is continuous functionals on the \Fre space $\Sc(U)$
of Schwartz functions on $U$. In this case, $\Sc(U)$ can be
identified with $\Sc(\R)$ via the stereographic projection.

The space of Schwartz functions on $\R$ is defined to be the space
of all infinitely differentiable functions that rapidly decay at
infinity together with all their derivatives, i.e. $x^nf^{(k)}$ is
bounded for any $n,k$.

The goal of this paper is to extend the notions of Schwartz
functions and Schwartz distributions to a larger geometric realm.

%\subsection{}
As we can see, the definition is of algebraic nature. Hence it
would not be reasonable to try to extend it to arbitrary smooth
manifolds. However, it is reasonable to extend this notion to
smooth algebraic varieties. Unfortunately, sometimes this is not
enough. For example, a connected component of real algebraic
variety is not always an algebraic variety. By this reason we
extend this notion to smooth semi-algebraic manifolds. They are
called \textbf{Nash manifolds}\footnote{The necessary
preliminaries on Nash manifolds are given in sections
\ref{SemiAlgGeo} and \ref{Nash}.}.

For any Nash manifold $M$, we will define the spaces $\G(M)$,
$\T(M)$ and $\Sc(M)$ of \textbf{generalized Schwartz
functions}\footnote{In this paper we distinguish between the
(similar) notions of a generalized function and a distribution.
They can be identified by choosing a measure. Without fixing a
measure, a smooth function defines a generalized function but not
a distribution. We will discuss it later in more details.},
\textbf{tempered functions} and \textbf{Schwartz functions} on
$M$. Informally, $\T(M)$ is the ring of functions that have no
more than polynomial growth together with all their derivatives,
$\G(M)$ is the space of generalized functions with no more than
polynomial growth and $\Sc(M)$ is the space of functions that
decay together with all their derivatives faster than any inverse
power of a polynomial.

As in the classical case, in order to define generalized Schwartz
functions, we have to define Schwartz functions first. Both
$\G(M)$ and $\Sc(M)$ are modules over $\T(M)$.

The triple $\Sc(M)$, $\T(M)$, $\G(M)$ is analogous to
$C_c^{\infty}(M)$, $C^{\infty}(M)$ and $C^{-\infty}(M)$ but it has
additional nice properties as we will see later.

We will show that for $M=\R^n$, $\Sc(M)$ is the space of classical
Schwartz functions and $\G(M)$ is the space of classical
generalized Schwartz functions. For compact Nash manifold $M$,
$\Sc(M)=\T(M)=C^{\infty}(M)$.

\subsection{Main results}

In this subsection we summarize the main results of the paper.
\begin{result}
Let $M$ be a Nash manifold and $Z\subset M$ be a closed Nash
submanifold. Then the restriction maps $\T(M) \rightarrow \T(Z)$
and $\Sc(M) \rightarrow \Sc(Z)$ are onto (see theorems \ref{ext}
and \ref{ExtTemp}).
\end{result}

\begin{result} \label{mainresult}
Let $M$ be a Nash manifold and $U \subset M$  be a semi-algebraic
open subset. Then a Schwartz function on $U$ is the same as a
Schwartz function on $M$ which vanishes with all its derivatives
on $M \setminus U$ (see theorem \ref{openset}).
\end{result}

This theorem tells us that extension by zero $\Sc(U)\rightarrow
\Sc(M)$ is a closed imbedding, and hence restriction morphism
$\G(M)\rightarrow \G(U)$ is onto.
\\Classical generalized functions do not have this property. This was our
main motivation for extending the definition of Schwartz
functions.

\subsection{Schwartz sections of Nash bundles}

Similar notions will be defined for Nash bundles, i.e. smooth
semi-algebraic bundles.

For any Nash bundle $E$ over $M$ we will define the spaces
$\G(M,E)$, $\T(M,E)$ and $\Sc(M,E)$ of generalized Schwartz,
tempered and Schwartz sections of $E$.

As in the classical case, a generalized Schwartz function is not
exactly a functional on the space of Schwartz functions, but a
functional on Schwartz densities, i.e. Schwartz sections of the
bundle of densities.

Therefore, we will define generalized Schwartz sections by
$\G(M,E)=(\Sc(M,\widetilde {E}))^*$, where $\widetilde
{E}=E^*\otimes D_M$ and $D_M$ is the bundle of densities on M.

Let $Z \subset M$ be a closed Nash submanifold, and $U=M \setminus
Z$. Result \ref{mainresult} tells us that the quotient space of
$\G(M)$ by the subspace $\G(M)_Z$ of generalized Schwartz
functions supported in $Z$ is $\G(U)$. Hence it is useful to study
the space $\G(M)_Z$. As in the classical case, $\G(M)_Z$ has a
filtration by the degree of transversal derivatives of delta
functions. The quotients of the filtration are generalized
Schwartz sections over $Z$ of symmetric powers of normal bundle to
$Z$ in $M$, after a twist.

This result can be extended to generalized Schwartz sections of
arbitrary Nash bundles (see corollary \ref{c4}).

\subsection{Restricted topology and sheaf properties}

Similarly to algebraic geometry, the reasonable topology on Nash
manifolds to consider is a topology in which open sets are open
semi-algebraic sets. Unfortunately, it is not a topology in the
usual sense of the word, it is only what is called
\textbf{restricted topology}. This means that the union of an
infinite number of open sets does not have to be open. The only
open covers considered in the restricted topology are finite open
covers.

The restriction of a generalized Schwartz function (respectively
of a tempered function) to an open subset is again a generalized
Schwartz (respectively a tempered function). This means that they
form pre-sheaves. We will show that they are actually sheaves,
which means that for any finite open cover $M = \bigcup \limits
_{i=1}^n U_i$, a function $\alpha$ on $M$ is tempered if and only
if $\alpha|_{U_i}$ is tempered for all $i$. It is of course not
true for infinite covers. For definitions of a pre-sheaf and a
sheaf in the restricted topology see section \ref{ResTop}. We
denote the sheaf of generalized Schwartz functions by $\G_M$ and
the sheaf of tempered functions by $\T_M$. By result
\ref{mainresult}, $\G_M$ is a flabby sheaf.

Similarly, for any Nash bundle $E$ over $M$ we will define the
sheaf $\T_M^E $ of tempered sections and the sheaf $\G_M^E$ of
generalized Schwartz sections.

As we have mentioned before, Schwartz functions behave similarly
to compactly supported smooth functions. In particular, they
cannot be restricted to an open subset, but can be extended by
zero from an open subset. This means that they do not form a
sheaf, but an object dual to a sheaf, a so-called cosheaf. The
exact definition of a cosheaf will be given in the appendix
(section \ref{cosheaf}). We denote the cosheaf of Schwartz
functions by $\Sc_M$.
%we will show that if $U$ is an open (semi-algebraic) subset of
%$M$, the extension by zero of a Schwartz function on $U$ is a
%Schwartz function on $M$.
%This property shows us the structure of
%cosheaf on Schwartz functions on $M$. Cosheaf is a notion dual to
%the notion of sheaf, i.e. instead of restriction to open subset we
%have extension from open subset.
We will prove that $\Sc_M$ is actually a cosheaf and not just
pre-cosheaf by proving a Schwartz version of the partition of
unity theorem. Similarly, for any Nash bundle $E$ over $M$ we will
define the cosheaf $\Sc_M^E$ of Schwartz sections.

\subsection{Possible applications}

Schwartz functions are used in the representation theory of
algebraic groups. Our definition coincides with Casselman's
definition (cf. \cite{Cas}) for algebraic groups. Our paper allows
to use Schwartz functions in additional situations in the
representation theory of algebraic groups, since an orbit of an
algebraic action is a Nash manifold, but does not have to be an
algebraic group or even an algebraic variety.

Generalized Schwartz sections can be used for ``devisage''. We
mean the following. Let $U\subset M$ be an open (semi-algebraic)
subset. Instead of dealing with generalized Schwartz sections of a
bundle on $M$, we can deal with generalized Schwartz sections of
its restriction to $U$ and generalized Schwartz sections of some
other bundles on $M \setminus U$ (see \ref{c4}).

For example if we are given an action of an algebraic group $G$ on
an algebraic variety $M$, and a $G$-equivariant bundle $E$ over
$M$, then devisage to orbits helps us to investigate the space of
$G$-invariant generalized sections of $E$. One of the
implementations of this argument appears in \cite{AGS}. There we
also use the fact that $\Sc(\R^n)$ is preserved by Fourier
transform.

\subsection{Summary} \label{Summary}

To sum up, for any Nash manifold $M$ we define a sheaf $\T_M$ of
algebras on $M$ (in the restricted topology) consisting of
tempered functions, a sheaf $\G_M$ of modules over $\T_M$
consisting of generalized Schwartz functions, and a cosheaf
$\Sc_M$ of modules over $\T_M$ consisting of Schwartz functions.

Moreover, for any Nash bundle $E$ over $M$ we define sheaves
$\T^E_M$ and $\G^E_M$ of modules over $\T_M$ consisting of
tempered and generalized Schwartz sections of $E$ respectively
and a cosheaf $\Sc^E_M$ of modules over $\T_M$ consisting of
Schwartz sections of $E$.

Let us list the main properties of these objects that we will
prove in this paper:

\begin{property}\label{p1}
$\!.$ Compatibility: For open semi-algebraic subset $U \subset M$,
$\Sc^E_M|_U=\Sc^{E|_U}_U$, $\T^E_M|_U=\T^{E|_U}_U$,
$\G^E_M|_U=\G^{E|_U}_U$.
\end{property}
\begin{property}\label{p2}
$\!.$ $\Sc(\R ^n)$ = Classical Schwartz functions on $\R ^n$ (see
corollary \ref{Usual}).
\end{property}
\begin{property}\label{p3}
$\!.$ For compact $M$, $\Sc(M,E) = \T(M,E)=$ $C^{\infty}(M,E)$
(see theorem \ref{ProofProp3}).
\end{property}
\begin{property}\label{p4}
$\!.$ $\G_M^{E}=(\Sc^{\widetilde {E} }_M)^*$ , where $\widetilde
{E}=E^*\otimes D_M$ and $D_M$ is the bundle of densities on M .
\end{property}
\begin{property}\label{p5}
$\!.$ Let $Z \subset M$ be a closed Nash submanifold. Then the
restriction maps  $\Sc(M,E)$ onto $\Sc(Z,E|_Z)$ and $\T(M,E)$ onto
$\T(Z,E|_Z)$ (see section \ref{BasProp}).
\end{property}
\begin{property}\label{p6}
$\!.$ Let $U \subset M$  be a semi-algebraic open subset, then
$$\Sc^E_M(U) \cong \{\phi \in \Sc(M,E)| \quad \phi \text{ is 0 on } M
\setminus U \text{ with all derivatives} \}.$$ (see theorem
\ref{openset}).
\end{property}
\begin{property} \label{p7}
$\!.$ Let $Z \subset M$ be a closed Nash submanifold. Consider
$\G(M,E)_Z=\{\xi \in \G(M,E) |\xi$ is supported in $Z \}$. It has
a canonical filtration such that its factors are canonically
isomorphic to $\G(Z,{E|_Z \otimes S^i(N_Z^M) \otimes {D_M^*|}_Z
\otimes D_Z})$ where $N_Z^M$ is the normal bundle of $Z$ in $M$
and $S^i$ means i-th symmetric power (see corollary \ref{c4}).
\end{property}
\subsection{Remarks}
\begin{remark}
Harish-Chandra has defined a Schwartz space for every reductive
Lie group. However, Harish-Chandra's Schwartz space does not
coincide with the space of Schwartz functions that we define in
this paper even for the algebraic group $\R^{\times}$.

%Harish-Chandra gave a different definitions of Schwartz functions,
%generalized Schwartz functions and tempered functions on Lie
%groups (cf. \cite{HC}).
%
%These definitions are not equivalent to our definition, even for
%the algebraic group $\R^{\times}$.
%Our definition does not coincide with Harish-Chandra's definition
%(cf. \cite{HC}) for Lie groups, because a Lie group isomorphism
%which is not a Nash diffeomorphism maps Harish-Chandra's Schwartz
%functions to Harish-Chandra's Schwartz functions but will not map
%our Schwartz functions to our Schwartz functions.
\end{remark}
\begin{remark}
There is a different approach to the concept of Schwartz
functions. Namely, if $M$ is embedded as an open subset in a
compact manifold $K$ then one can define the space of Schwartz
functions on $M$ to be the space of all smooth functions on $K$
that vanish outside $M$ together with all their derivatives. This
approach is implemented in \cite{CHM}, \cite{KS}, \cite{Mor} and
\cite{Pre}. In general, this definition depends on the embedding
into $K$. Our results show that for Nash manifolds $M$ and $K$ it
coincides with our definition and hence does not depend on the
embedding.
\end{remark}

\begin{remark}
After the completion of this project we found out that many of the
properties of Schwartz functions on affine Nash manifolds, that is
most of section \ref{AffSchwartz} have been obtained already in
\cite{dCl}.
\end{remark}

\begin{remark}
In fact, tempered functions can be defined in terms of Schwartz
functions or in terms of generalized Schwartz functions by
$$\T(M,E) = \{\alpha \in
C^{\infty}(M,E)|\,\,\, \alpha \Sc(M)\subset \Sc(M,E)\} =$$
$$\quad \quad \quad =\{\alpha \in C^{\infty}(M,E) |\,\,\, \alpha \G(M)\subset
\G(M,E)\},$$ but the proof is rather technical and lies out of the
scope of this paper.
\end{remark}
\begin{remark}
Throughout the whole paper ``smooth'' means infinitely
differentiable and ``bounded'' means bounded in absolute value,
unless stated otherwise.
\end{remark}
\subsection*{Acknowledgements}

We would like to thank our teacher { Prof. Joseph Bernstein}
for teaching us most of the mathematics we know and for his help
in this work.

We would also like to thank { Prof. William Casselman} who
dealt with similar things in the early 90s and sent us his
unpublished material.

We thank { Prof. Michel Duflo} for drawing our attention to the work
\cite{dCl} and { Prof. Pierre Schapira} for drawing our attention to
the works \cite{KS}, \cite{Mor} and \cite{Pre} which are related
to this paper.

We would like to thank  { Prof. Semyon Alesker}, { Lev
Buhovski}, { Prof. Alexander Braverman}, { Vadim Kosoy},
{ Prof. Pierre Milman}, { Prof. Vitali Milman}, { Dr.
Omer Offen}, { Dr. Eitan Sayag}, { Dr. Michael Temkin} and
{ Dr. Ilya Tyomkin} for helpful discussions.

Also we would like to thank { Dr. Ben-Zion Aizenbud}, {
Shifra Reif}, {Ilya Surdin} and { Frol Zapolsky} for proof
reading.

Finally, we would like to thank {Prof. Paul Biran}, { Prof.
Vladimir Berkovich}, { Prof. Stephen Gelbart},  { Prof.
Eugenii Shustin},  { Prof. Sergei Yakovenko} for their useful
remarks.
\section{Semi-algebraic geometry} \label{SemiAlgGeo}

In this section we give some preliminaries on semi-algebraic
geometry from \cite{BCR}.
%, \cite{DK} and \cite{Shi} and some of
%their technical extensions that will be necessary for us to
%proceed.
%
\subsection{Basic notions}
\begin{definition}
{A subset $A \subset \R^n$ is called a \textbf{semi-algebraic set}
iff it can be presented as a finite union of sets defined by a
finite number of polynomial equalities and inequalities. In other
words, there exist finitely many polynomials $f_{ij}, g_{ik} \in
\R[x_1,\dots,x_n]$ such that $$A = \bigcup \limits _{i=1}^r \{x
\in \R^n | f_{i1}(x)>0,\dots,f_{is_i}(x)>0,
g_{i1}(x)=0,\dots,g_{it_i}(x)=0\}.$$ }
\end{definition}
\begin{lemma}
The collection of semi-algebraic sets is closed with respect to
finite unions, finite intersections, and complements.
\end{lemma}
The proof is immediate.
\begin{definition}
{Let $A \subset \R^n $ and $B \subset \R^m$ be semi-algebraic
sets. A mapping $\nu:A \rightarrow B $ is called
\textbf{semi-algebraic} iff its graph is a semi-algebraic subset
of $\R^{m+n}$.}
\end{definition}
\subsection{Tarski-Seidenberg principle of quantifier elimination \\ and its applications}
One of the main tools in the theory of semi-algebraic spaces is
the Tarski-Seidenberg principle of quantifier elimination. Here we
will formulate and use a special case of it. We start from the
geometric formulation.
\begin{theorem}
Let $A \subset \R^n$ be a semi-algebraic subset and $p:\R^n
\rightarrow \R^{n-1}$ be the standard projection. Then the image
$p(A)$ is a semi-algebraic subset of $\R^{n-1}$.
\end{theorem}
By induction and a standard graph argument we get the following
corollary.
\begin{corollary}
An image of a semi-algebraic subset of $\R^n$ under a
semi-algebraic map is semi-algebraic.
\end{corollary}
Sometimes it is more convenient to use the logical formulation of
the Tarski-Seidenberg principle. Informally it says that any set
that can be described in semi-algebraic language is
semi-algebraic. We will now give the logical formulation and
immediately after that define the logical notion used in it.
\begin{theorem} {(Tarski-Seidenberg principle)}
Let $\Phi$ be a formula of the language $L(\R)$ of ordered fields
with parameters in $\R$. Then there exists a quantifier - free
formula $\Psi$ of $L(\R)$ with the same free variables
$x_1,\dots,x_n$ as $\Phi$ such that
 $\forall x \in \R^n, \Phi(x) \Leftrightarrow \Psi(x)$.
\end{theorem}
For the proof see Proposition 2.2.4 on page 28 of \cite{BCR}.
\begin{definition}
{A \textbf{first-order formula of the language of ordered fields
with parameters in $\R$} is a formula written with a finite number
of conjunctions, disjunctions, negations and universal and
existential quantifiers ($\forall$ and $\exists$) on variables,
starting from atomic formulas which are formulas of the kind
$f(x_1,\dots,x_n) = 0$ or $g(x_1,\dots,x_n) > 0$, where $f$ and
$g$ are polynomials with coefficients in $\R$. The free variables
of a formula are those variables of the polynomials which are not
quantified. We denote the language of such formulas by $L(\R)$.}
\end{definition}
\begin{notation}
Let $\Phi$ be a formula of $L(\R)$ with free variables
$x_1,\dots,x_n$. It defines the set of all points
$(x_1,\dots,x_n)$ in $\R^n$ that satisfy $\Phi$. We denote this
set by $\mathbf{S_{\Phi}}$. In short, $$\mathbf{S_{\Phi}}:=\{ x
\in \R^n |\Phi(x)\}.$$ %The next corollary says that it is semi-algebraic.
\end{notation}
\begin{corollary}
Let $\Phi$ be a formula of $L(\R)$. Then $S_{\Phi}$ is a
semi-algebraic set.
\end{corollary}
\emph{Proof.} Let $\Psi$ be a quantifier-free formula equivalent
to $\Phi$. The set $S_{\Psi}$ is semi-algebraic since it is a
finite union of sets defined by  polynomial equalities and
inequalities. Hence $S_{\Phi}$ is also semi-algebraic since
$S_{\Phi}=S_{\Psi}$. \proofend
\begin{proposition}
The logical formulation of the Tarski-Seidenberg principle implies
the geometric one.
\end{proposition}
\emph{Proof.} Let $A \subset \R^n$ be a semi-algebraic subset, and
$pr:\R^n \rightarrow \R^{n-1}$ the standard projection. Then there
exists a formula $\Phi \in L(\R)$ such that $A = S_{\Phi}$. Then
$pr(A)= S_{\Psi}$ where
$$\Psi(y)= \text{``}\exists x \in \R^n \,(\pr(x) = y \wedge \Phi(x))".$$ Since
$\Psi \in L(\R)$, the claim follows from the previous corollary.
\begin{remark}
In fact, it is not difficult to deduce the logical formulation
from the geometric one.
\end{remark}
Let us now demonstrate how to use the logical formulation of the
Tarski-Seidenberg principle.
\begin{corollary}
The closure of a semi-algebraic set is semi-algebraic.
\end{corollary}
\emph{Proof.} Let $A \subset \R^n$ be a semi-algebraic subset, and
let $\overline{A}$ be its closure. Then $\overline{A}=S_{\Psi}$
where
$$\Psi(x)= "\forall \eps>0 \, \exists y \in A \, |x-y|^2<\eps".$$
Clearly, $\Psi \in L(\R)$ and hence $\overline{A}$ is
semi-algebraic. \proofend
\begin{corollary}
Images and preimages of semi-algebraic sets under semi-algebraic
mappings are semi-algebraic.
\end{corollary}
%\emph{Proof.} This proposition follows from the Tarski-Seidenberg
%principle. \proofend
%
\begin{proposition}
Let $\nu$ be a bijective semi-algebraic mapping. Then the inverse
mapping $\nu^{-1}$ is also semi-algebraic.
\end{proposition}
\emph{Proof.} The graph of $\nu$ is obtained from the graph of
$\nu^{-1}$ by switching the coordinates. $\text{ }$ \proofend
\begin{proposition}
$ $\\(i) The composition of semi-algebraic mappings is
semi-algebraic.\\
(ii) The $\R$-valued semi-algebraic functions on a semi-algebraic
set $A$ form a ring, and any nowhere vanishing semi-algebraic
function is invertible in this ring.
\end{proposition}
\emph{Proof.} \\
(i) Let $\mu: A \rightarrow B$ and $\nu: B \rightarrow C$ be
semi-algebraic mappings. Let $\Gamma_{\mu} \subset \R^{m+n}$ be
the graph of $\mu$ and $\Gamma_{\nu} \subset \R ^{n+p}$ be the
graph of $\nu$. The graph of $\nu \circ \mu$ is the projection of
$(\Gamma_{\mu} \times \R^p)\cap (\R ^m \times
\Gamma_{\nu})$ onto $\R ^{m+p}$ and hence is semi-algebraic.\\
(ii) follows from (i) by noting that $F+G$ is the composition of
$(F,G):A \rightarrow \R^2$ with $+:\R^2 \rightarrow \R$, $FG$ is
the composition of $(F,G):A \rightarrow \R^2$ with $\times:\R^2
\rightarrow \R$ and $\frac{1}{F}$ is the composition of $F$ with
$\frac{1}{x}: \R \setminus 0 \rightarrow \R \setminus 0 $.
\proofend
\begin{lemma} \label{majalg}
Let $X\subset \R^n$ be a closed semi-algebraic subset. Then any
continuous semi-algebraic function $F:X \rightarrow \R$ can be
majorated \footnote{By majorated we mean bounded by absolute value
from above.} by the restriction to $X$ of some polynomial on
$\R^n$.
\end{lemma}
\emph{Proof} Let $G:\R \rightarrow \R$ be defined by $G(r) := \max
\limits _{\{x\in X| |x| \leq r\}}F(x)$. By Tarski-Seidenberg
principle, $G$ is semi-algebraic. The lemma now reduces to its
one-dimensional case, which is proven on page 43 of \cite{BCR}
(proposition 2.6.2).
 \proofend
\subsection{Additional preliminary results}
%Now we will formulate a couple of theorems that will be useful for
%us.
\begin{theorem} \label{saloc}
Let $F:A \rightarrow \R$ be a semi-algebraic function on a locally
closed semi-algebraic set. Let $Z(F):= \{x\in A | F(x) = 0 \}$ be
the set of zeros of $F$ and let $A_F:=A\setminus Z(F)$ be its
complement. Let $G:A_F \rightarrow \R$ be a semi-algebraic
function. Suppose that $F$ and $G$ are continuous. Then there
exists an integer $N>0$ such that the function $F^NG$, extended by
0 to $Z(F)$, is continuous on $A$.
\end{theorem}
The proof can be found on page 43 of \cite{BCR} (proposition
2.6.4).
\begin{theorem}{(Finiteness).} \label{finbas}
Let $X \subset \R^n$ be a semi-algebraic set. Then every open
semi-algebraic subset of $X$ can be presented as a finite union of
sets of the form $\{x\in X |\, p_i(x)>0, i=1 \dots n\}$, where
$p_i$ are polynomials in $n$ variables.
\end{theorem}
The proof can be found on page 46 of \cite{BCR} (theorem 2.7.2).
\begin{theorem} \label{Components}
Every semi-algebraic subset of $\R^n$ has a finite number of
connected components, which are semi-algebraic
\end{theorem}
The proof can be found on page 35 of  \cite {BCR} (theorem 2.4.5).
\section{Nash manifolds} \label{Nash}
In this section we give some preliminaries on Nash
manifolds\footnote{What we mean by Nash manifold is sometimes
called $C^{\infty}$ Nash manifold or $C^{\omega}$ Nash manifold.}
from \cite{BCR}, \cite{DK} and \cite{Shi} and some of their
technical extensions that will be necessary for us to proceed.
%Now we would like to define the category of Nash manifolds.
Most of section \ref{AffSchwartz} does not rely on subsections
\ref{ResTop}-\ref{NashTubNeigh}. However, section \ref{SchAbsNash}
will use all of this section.

The theory of Nash manifolds is similar both to differential
topology and algebraic geometry. Our approach to Nash manifolds
comes from algebraic geometry.

%In subsection \ref{NashSubman} we will define

\subsection{Nash submanifolds of $\R^n$} \label{NashSubman}
\begin{definition}
{A \textbf{Nash map} from an open semi-algebraic subset $U$ of
$\R^n$ to an open semi-algebraic subset $V$ of $\R^m$ is a smooth
(i.e. infinitely differentiable) semi-algebraic map. The ring of
$\R$-valued Nash functions on $U$ is denoted by ${\cal N}(U)$. A
\textbf{Nash diffeomorphism} is a Nash bijective map whose inverse
map is also Nash. }
\end{definition}
\begin{remark}
In fact, a Nash map is always real analytic (cf. \cite{Mal} or
\cite{Shi}, Corollary I.5.7) but we will not use this.
\end{remark}
\begin{remark}
A Nash map which is a diffeomorphism is a Nash diffeomorphism,
since the inverse of a semi-algebraic map is semi-algebraic. Note
also that a partial derivative of Nash function is Nash by the
Tarski-Seidenberg principle.
\end{remark}
As we are going to work with semi-algebraic differential geometry,
we will need a semi-algebraic version of implicit function
theorem.
\begin{theorem}(Implicit Function Theorem) Let $(x^0,y^0) \in \R^{n+p}$ and let $f_1,\dots,f_p$ be Nash
functions on an open neighborhood of $(x^0,y^0)$ such that
$f_j(x^0,y^0)=0$ for $j=1,..,p$ and the matrix $[\frac{\partial
f_j}{\partial y_i}(x^0,y^0)]$ is invertible. Then there exist open
semi-algebraic neighborhoods $U$ and $V$ of $x^0$ (resp. $y^0$) in
$\R^n$ (resp. $\R^p$) and a Nash mapping $\nu$, such that
$\nu(x^0)=y^0$ and $f_1(x,y) = \dots = f_p(x,y)=0 \Leftrightarrow
y = \nu(x)$ for every $(x,y) \in U \times V.$
\end{theorem}
The proof can be found on page 57 of \cite{BCR} (corollary 2.9.8).
\begin{definition}
{ A \textbf{Nash submanifold of} $\mathbf{\R^n}$ is a
semi-algebraic subset of $\R^n$ which is a smooth submanifold.}
\end{definition}
By the implicit function theorem it is easy to see that this
definition is equivalent to the following one, given in
\cite{BCR}:

\begin{definition}
{A semi-algebraic subset $M$ of $\R^n$ is said to be a
\textbf{Nash submanifold of $\mathbf{\R^n}$ of dimension $d$} if,
for every point $m$ of $M$, there exists a Nash diffeomorphism
$\nu$ from an open semi-algebraic neighborhood $\Omega$ of the
origin in $\R^n$ onto an open semi-algebraic neighborhood
$\Omega'$ of $m$ in $\R^n$ such that $\nu(0) = m$ and $\nu(\R^d
\times \{0\}) \cap \Omega) = M \cap \Omega'$. }
\end{definition}
\begin{theorem}\footnote{ We will use this theorem since it makes
several formulations and proofs shorter, but in each case it can
be avoided.} \label{OpenAffine} Any Nash submanifold of $\R^n$ is
Nash diffeomorphic to a closed Nash submanifold of $\R^N$.
\end{theorem}

For proof see Corollary I.4.3 in \cite{Shi} or theorems 8.4.6 and
2.4.5 in \cite{BCR}.
\begin{definition}
{A \textbf{Nash function} on a Nash submanifold $M$ of $\R^n$ is a
semi-algebraic smooth function on $M$. The ring of $\R$-valued
Nash functions on $M$ is denoted by ${\cal N}(M)$. }
\end{definition}

The rest of section \ref{Nash} is not necessary for readers
interested only in Schwartz functions on affine Nash manifolds.
\subsection{Restricted topological spaces and sheaf theory over them.} \label{ResTop}
Now we would like to define Nash manifolds independently of their
embedding into $\R^n$. Analogously to algebraic geometry we will
define them as ringed spaces. Hence we will need to introduce
topology and structure sheaf on Nash manifolds. The natural
topology to consider is topology of open semi-algebraic sets.
Unfortunately, infinite unions of semi-algebraic sets are not
necessary semi-algebraic, hence it is not a topology in the usual
sense of the word. Therefore, we will need to define a different
notion of topology and introduce an appropriate sheaf theory over
it. In this section we follow \cite{DK}. A similar use of
restricted topology appears in \cite{Mor} and \cite{Pre}.
\begin{definition}
{A \textbf{restricted topological space $M$} is a set $M$ equipped
with a family $ \overset {_{\circ}} {\mathfrak S}   (M)$ of
subsets of $M$, called the open subsets which contains $M$ and the
empty set, and is closed with respect to \textbf{finite} unions
and finite intersections.}
\end{definition}
\begin{remark}
Pay attention that in general, there is no closure in restricted
topology since infinite intersection of closed sets does not have
to be closed. In our case, open sets will have closure.
\end{remark}
\begin{remark}
A restricted topological space $M$ can be considered as a site in
the sense of Grothendieck. The category of the site has the open
sets of $M$ as objects and the inclusions as morphisms. The covers
$(U_i \rightarrow U)_{i\in I}$ are the \textbf{finite} systems of
inclusions with $\bigcup \limits _{i=1}^n U_i = U$. The standard
machinery of Grothendieck topology gives us the notions of a
pre-sheaf and a sheaf on $M$. Now we will repeat the definitions
of these notions in simpler terms.
\end{remark}
\begin{definition}
A \textbf{pre-sheaf $F$} on a restricted topological space $M$ is
an assignment $U \mapsto F(U)$ for every open $U$ of an abelian
group, vector space, etc., and for every inclusion of open sets $V
\subset U$ a restriction morphism $res_{U,V}: F(U) \rightarrow
F(V)$ such that $res_{U,U} = Id$ and for $W \subset V \subset U$,
$res_{V,W} \circ res_{U,V} = res_{U,W}$.
\end{definition}
\begin{definition}
{A \textbf{sheaf $F$} on a restricted topological space $M$  is a
pre-sheaf fulfilling the usual sheaf conditions, except that now
only finite open covers are admitted. The conditions are: for any
open set $U$ and any finite cover $U_i$ of $M$ by open subsets,
the sequence $$0 \rightarrow F(U) \overset{res_1}{\rightarrow}
\prod_{i=1}^n F(U_i) \overset{res_2}{\rightarrow}
\prod_{i=1}^{n-1} \prod_{j=i+1}^{n} F(U_i \cap U_j)$$ is exact.\\
The map $res_1$ above is defined by $res_1(\xi) = \prod \limits
_{i=1}^n res_{U,U_i}(\xi)$ and the map $res_2$ by
$$res_2(\prod_{i=1}^n \xi_i) = \prod_{i=1}^{n-1} \prod_{j=i+1}^n res_{U_i,U_i \cap U_j}(\xi_i)
- res_{U_j,U_i \cap U_j}(\xi_j)$$ . }
\end{definition}
%\newpage
\begin{definition}
An \textbf{$\R$-space} is a pair $(M,{\cal O}_M)$ where $M$ is a
restricted topological space and ${\cal O}_M$ a sheaf of
$\R$-algebras over $M$ which is a subsheaf of the sheaf $C_M$ of
all continuous real-valued functions on $M$.

A \textbf{morphism between $\R$-spaces} $(M,{\cal O}_M)$ and
$(N,{\cal O}_N)$ is a continuous map $\nu:M \rightarrow N$, such
that for every open set $U\subset N$ and every function $f \in
{\cal O}_N(U)$, the composition $f \circ \nu|_{\nu^{-1}(U)}$ lies
in ${\cal O}_M(\nu^{-1}(U))$.
%the induced morphism of sheaves $\nu^*:\nu^*(\R\{N\}) ??? \nu^* twise - difficult to understand
%\rightarrow \R\{M\}$ maps ${\cal O}_N$ to ${\cal O}_M$.}
\end{definition}

\begin{remark}
In a dual way, one can define the notion of a cosheaf over a
restricted topological space. We will use this notion only in
section \ref{SchAbsNash} and give the exact definition in
subsection \ref{cosheaf}.
\end{remark}

\subsection{Abstract Nash manifolds}

\begin{example}
To any Nash submanifold $M$ of $\R^n$ we associate an $\R$-space
in the following way. Take for $\overset {_{\circ}} {\mathfrak S}
(M)$ the family of all open subsets of $M$ which are
semi-algebraic in $\R^n$. For any open (semi-algebraic) subset $U$
of $M$ we set ${\cal O}_M(U)$ to be the algebra ${\cal N}(U) $ of
Nash functions $U \rightarrow \R$.
\end{example}
\begin{definition}
{An \textbf{affine Nash manifold} is an $\R$-space which is
isomorphic to an $\R$-space associated to a closed Nash
submanifold of $\R^n$.}
\end{definition}
\begin{definition}
{A \textbf{Nash manifold} is an $\R$-space $(M,{\cal N}_M)$ which
has a finite cover $(M_i)_{i=1}^n$ by open sets $M_i$ such that
the $\R$-spaces $(M_i,{\cal N}_M|_{M_i})$ are affine Nash
manifolds. A morphism between Nash manifolds is a morphism of
$\R$-spaces between them. }
\end{definition}
\begin{remark}
A map between two closed Nash submanifolds of $\R^n$ is Nash if
and only if it is a morphism of Nash manifolds. Hence we will call
morphisms of Nash manifolds Nash maps, and isomorphisms of Nash
manifolds Nash diffeomorphisms.
\end{remark}
The following proposition is a direct corollary of theorem
\ref{OpenAffine}.
\begin{proposition} \label{OpenNash}
1. Any open (semi-algebraic) subset $U$ of an affine Nash manifold
$M$ is an affine Nash manifold.\\
2. Any open (semi-algebraic) subset $U$ of a Nash manifold $M$ is
a Nash manifold.
\end{proposition}
%\emph{Proof.} It is enough to prove for the case when $M$ is
%affine, i.e. $M$ is a closed Nash submanifold of $\R^n$. By
%finiteness theorem (\ref{finbas}), it is enough to prove the
%proposition for $U=\{x \in M | f_i(x)>0$ for $1 \leq i \leq k \}$,
%where $f_i$ are polynomials on $\R^n$. For this case it is done
%the same way as it is done in algebraic geometry. Consider the set
%$M_{f_k}:= \{(x_1,\dots,x_{n+1}) \in \R^{n+1} | (x_1,\dots,x_n) \in M$
%and $f_k(x_1,\dots,x_{n})x_{n+1}=1$ and $x_{n+1} \geq 0 \}$. It is
%clearly an affine Nash variety. As a ringed space, $U$ can be
%naturally embedded into $M_{f_k}$ and defined there by $k-1$
%strict inequalities. Now, by induction on $k$, $U$ is a Nash
%manifold. \proofend
%
\subsubsection{Examples and Remarks}
\begin{example}
Any smooth affine algebraic variety over $\R$ is an affine Nash
manifold.
\end{example}
\begin{remark}
A union of several connected components of an affine Nash manifold
is an affine Nash manifold. This is true since every
semi-algebraic subset of $\R^n$ has a finite number of connected
components and each of them is semi-algebraic (see \cite{BCR},
theorem 2.4.5).
\end{remark}
\begin{remark}
Any affine Nash manifold is Nash diffeomorphic to a smooth real
affine algebraic variety (see Chapter 14 of \cite{BCR} for compact
Nash manifold, and remark VI.2.11 in \cite{Shi} for non-compact
Nash manifold). However, the category of affine Nash manifolds is
richer than the category of smooth real affine algebraic
varieties, because it has more morphisms. In particular, two
non-isomorphic smooth real affine algebraic varieties can be Nash
diffeomorphic. For example, the hyperbola $\{xy=1\}$ is Nash
diffeomorphic to the union of two straight lines $\{x+y=1/2\}
\bigsqcup \{x+y=-1/2\}$. %which means $\R\setminus \{0\}$ and $\R \bigsqcup \R$.
\end{remark}
\begin{remark}
Note that the Nash groups $(\R_{> 0},\times)$ and $(\R,+)$ are
\textbf{not} isomorphic as Nash groups, although they are both
Nash diffeomorphic and isomorphic as Lie groups.
\end{remark}
\begin{remark}
Any quasiprojective Nash manifold is affine since any projective
Nash manifold is affine (see page 23 in \cite{Shi} after the proof
of lemma I.3.2).
\end{remark}
\begin{remark}
Any Nash manifold has an obvious natural structure of a smooth
manifold.
\end{remark}
%
%??? \subsubsection{}
%In this paper we will need some delicate work with covers. In
%differential topology there is a theorem that for any cover
%$M=\cup U_i$ there exists a cover $M=\cup V_j$ such that $V_j$ are
%basic open subsets (homeomorphic to $\R^n$) and $\forall j \exists
%i. \overline{V_j}\subset U_i$. we will prove a strong version of
%this theorem for Nash manifolds. For exact formulation and proof
%of this theorem see appendix \ref{App}, theorem \ref{UltNarCov}.
%In the appendix we will also give other technical results and
%definitions.

\subsection{Nash vector bundles}
\begin{definition}
{Let $\pi:M \rightarrow B$ be a Nash map of Nash manifolds. It is
called a \textbf{Nash locally trivial fibration} with fiber $Z$ if
$Z$ is a Nash manifold and there exist a \textbf{finite} cover $B
= \bigcup \limits _{i=1}^n U_i$ by open (semi-algebraic) sets and
Nash diffeomorphisms $\nu_i$ of $\pi^{-1}(U_i)$ with $U_i \times
Z$ such that the composition $\pi \circ \nu_i^{-1}$ is the natural
projection.}
\end{definition}
\begin{definition}
{Let $M$ be a Nash manifold. As in differential geometry, a
\textbf{Nash vector bundle $E$ over $M$} is a Nash locally trivial
fibration with linear fiber and such that trivialization maps
$\nu_i$ are fiberwise linear. By abuse of notation, we use the
same letters to denote bundles and their total spaces.}
\end{definition}
\begin{definition}
{Let $M$ be a Nash manifold and $E$ a Nash bundle over $M$. A
\textbf{Nash section of $E$} is a section of $E$ which is a Nash
map.}
\end{definition}
\begin{remark} \label{VeryNash}
In some books, for example \cite{BCR} such a bundle is called
pre-Nash. In these books, a bundle is called Nash if it can be
embedded into a trivial one. For affine manifold $M$ this property
implies that for open $U\subset M$ there exists a finite open
cover $U = \cup U_i$ such that for any $i$, any Nash section $s:
U_i \to E$ is a combination $s=\sum f_j t_j|_{U_i}$ where $t_j:M
\to E$ are global Nash sections and $f_j \in \mathcal{N}(U_i)$ are
Nash functions on $U_i$. For proof see \cite{BCR}.
\end{remark}
\begin{remark}
Direct sums, tensor products, external tensor products, the dual,
exterior powers, etc., of Nash vector bundles all have canonical
structures of Nash vector bundles.
\end{remark}
\begin{theorem}
Tangent, normal and conormal bundles, the bundle of differential
$k$-forms, etc., of a Nash manifold have canonical structures of
Nash bundles.
\end{theorem}
For the construction see section \ref{SecBunApp} in the Appendix.
\begin{theorem}
The bundle of densities of a Nash manifold $M$ has a canonical
structure of a Nash bundle.
\end{theorem}
For the construction see subsection \ref{BunDenApp} in the
Appendix.
\begin{notation}
We denote the bundle of densities of a Nash manifold $M$ by
\textbf{$D_M$}.
\end{notation}
\begin{notation}
{Let $E$ be a Nash bundle over $M$. We denote
\textbf{$\widetilde{E}$}$:=E^*\otimes D_M$.}
\end{notation}
\begin{definition}
{Let $M$ be a Nash manifold. Then a \textbf{Nash vector field on
$M$} is a Nash section of the tangent bundle of $M$. A
\textbf{Nash covector field on $M$} is a Nash section of the
cotangent bundle of $M$. A \textbf{Nash differential $k$-form on
$M$} is a Nash section of the bundle of differential $k$-forms on
$M$.}
\end{definition}
\begin{proposition} \label{remcot}
Let $M\subset \R^n$ be a closed Nash submanifold. Then the space
of Nash covector fields on $M$ is generated over ${\cal N}(M)$ by
$dx_i$.
\end{proposition}
For the proof see proposition \ref{BunAppCor} in the Appendix.
%
%
%
%\begin{definition}
%{\rm Let $f:M \rightarrow N$ be a Nash morphism of Nash manifolds.
%Let $\underset {M} { \overset {E}{\downarrow}}$ be a Nash bundle.
%Then $f^*E$ has a natural structure of a Nash bundle.}
%%We also
%%define \emph{twisted pull $f^?(E)$} by $f^?(E):= f^*(E \otimes
%%D_N*) \otimes D_M$.}
%\end{definition}
%%
%\begin{remark}
%Let $f:M \rightarrow N $ be a Nash diffeomorphism. Then
%differential defines canonical isomorphism $f^*(D_N) \cong D_M $.
%\end{remark}
%
\subsection{Nash differential operators}

As we have seen in the introduction, classical Schwartz functions
can be defined using polynomial differential operators. In order
to define Schwartz functions on Nash manifolds, we will use Nash
differential operators, that will be defined in this subsection.
\begin{definition}
Let $M$ be a smooth real affine algebraic variety. \textbf{The
algebra of algebraic differential operators on $M$} is the
subalgebra with 1 of $Hom_{\R}(C^\infty(M),C^\infty(M))$ generated
by multiplications by polynomial functions and by derivations
along algebraic vector fields.

Let $M$ be an affine Nash manifold. \textbf{The algebra of Nash
differential operators on $M$} is the subalgebra with 1 of
$Hom_{\R}(C^\infty(M),C^\infty(M))$ generated by multiplications
by Nash functions and by derivations
along Nash vector fields.%, where Nash vector field is a Nash section of the tangent bundle of $M$.
\end{definition}

The following lemma immediately follows from the Tarski-Seidenberg
theorem.
\begin{lemma} \label{DiffNashtoNash}
Any Nash differential operator maps ${\cal N}(M)$ to ${\cal
N}(M)$.
\end{lemma}

\begin{remark} \label{GrothDiffOp}
One could give equivalent definition of the algebra of Nash
differential operators on $M$ as Grothendieck's algebra of
differential operators over the algebra ${\cal N}(M)$.
\end{remark}

We will use the following trivial lemma

\begin{lemma} \label{DiffonOpen}
Let $U\subset \R^n$ be open (semi-algebraic) subset. Then any Nash
differential operator $D$ on $U$ can be written as $\sum_{i=1}^k
f_i(D_i|_U)$ where $f_i$ are Nash functions on $U$ and $D_i$ are
Nash differential operators on $\R^n$.
\end{lemma}

\begin{remark}
A similar but weaker lemma holds for any affine Nash manifold and
its open subset. It can be proven using remark \ref{VeryNash}.
\end{remark}
%
%\begin{lemma} \label{LocNashDiff}                                       %??? proof
%Let $M$ be an affine Nash manifold and $U$ be an affine open
%subset. Then the space of Nash differential operators on $U$ is
%generated as a module over ${\cal N}(U)$ by restrictions of Nash
%differential operators on $M$.
%\end{lemma}
%
\subsubsection{Algebraic differential operators on a Nash manifold}

The following ad hoc definition will be convenient for us for
technical reasons.

\begin{definition}
Given a closed embedding of $M$ to $\R^n$ we can define the notion
of \textbf{algebraic differential operators on $M$}, which depends
on the embedding. Algebraic vector fields on $M$ are defined as
Nash vector fields obtained by composition of the restriction to
$M$ and the orthogonal projection to $T_M$ from algebraic vector
fields on $\R^n$. The algebra of algebraic differential operators
on $M$ is the subalgebra with 1 of
$Hom_{\R}(C^\infty(M),C^\infty(M))$ generated by multiplications
by restriction of polynomial functions and by derivations along
algebraic vector fields.
\end{definition}

\begin{lemma} \label{multalg}
Let $M \subset \R^n$ be a closed affine Nash submanifold. Then the
space of Nash differential operators on $M$ is generated as a
module over ${\cal N}(M)$ by algebraic differential operators.
\end{lemma}
\emph{Proof.} Consider the algebraic vector fields
$\frac{\partial}{\partial x_i}|_M$ obtained from the standard
vector fields $\frac{\partial}{\partial x_i}$ on $\R^n$ by
restriction to $M$ and orthogonal projection to $T_M$. By
proposition \ref{remcot}, $\frac{\partial}{\partial x_i}|_M$
generate Nash vector fields over ${\cal N}(M)$. The lemma now
follows from lemma \ref{DiffNashtoNash}, the chain rule and the
Leibnitz rule. \proofend
%
%??? \begin{definition}
%{Let $f:M \rightarrow N$ be a Nash map of Nash manifolds. It
%is called \emph{projective} if it can be represented as a
%composition of a closed Nash embedding and a projection $\PP^n
%\times N \rightarrow N$}.
%\end{definition}
%
\subsection{Nash tubular neighborhood} \label{NashTubNeigh}

We will need a Nash analog of the tubular neighborhood theorem
from differential geometry.
\begin{notation}
{Let $M$ be a Nash manifold. Let $Z\subset M$ be a closed Nash
submanifold. We denote by $N_{Z}^{M}$ the normal bundle of $Z$ in
$M$. In case $M$ is equipped with a Nash Riemannian metric, we
consider $N_{Z}^{M}$ as a subbundle of the tangent bundle of $M$.}
\end{notation}
\begin{notation}
{Let $x \in \R^n$ and $r \in \R$. We denote the open ball with
center $x$ and radius $r$ by $B(x,r)$.

Let $M$ be a Nash manifold with a Nash Riemannian metric on it.
Let $Z\subset M$ be a closed Nash submanifold and $F:Z \rightarrow
\R$ a function on it. We define $B_M(Z,F)\subset N_Z^M$ by
$B_M(Z,F):= \{(z,v) \in N_Z^M|$ $||v||<F(z)\}$ }.
\end{notation}
\begin{remark}
Note that if $F$ is strictly positive and semi-algebraic then
$B_M(Z,F)$ is Nash.
\end{remark}
\begin{definition}
{Let $M$ be a Nash manifold with a Nash Riemannian metric on it.
Let $Z\subset M$ be a closed Nash submanifold. A \textbf{Nash
tubular neighborhood of $Z$ in $M$} is the following data: an open
Nash neighborhood $Tube(Z,M)$ of $Z$ in $M$, a strictly positive
Nash function $\rho_Z^M \in \mathcal{N}(Z)$ and a Nash
diffeomorphism $\nu_Z^M: B_M(Z,\rho_Z^M) \simeq Tube(Z,M)$.}
\end{definition}
\begin{theorem}{(Nash Tubular Neighborhood). } \label{NashTube}
Let $Z \subset M \subset \R^n$ be closed affine Nash submanifolds.
Equip $M$ with the Riemannian metric induced from $\R^n$. Then $Z$
has a Nash tubular neighborhood. Moreover, for any closed Nash
submanifold $Y\subset M$ disjoint from $Z$, $Z$ has a Nash tubular
neighborhood in $M\setminus Y$.
\end{theorem}

\emph{Proof.} Denote by $\nu_1$ the map $\nu_1:N_{Z}^{M}
\rightarrow \R^n$ given by $\nu_1(z,v) = z+v$. Denote by $W$ the
set of all points $x \in \R^n$ such that there exists a unique
point of $M$ closest to $x$. Denote $U := \overset{o}{W}$. It is
open semi-algebraic, and the projection $pr:U \rightarrow M$ is
also semi-algebraic. Denote $V = \nu_1^{-1}(U)$ and $\nu = pr
\circ \nu_1 :V \rightarrow M$.

From differential topology we know that there exists a smooth
strictly positive function $\sigma$ such that
$\nu|_{B_M(Z,\sigma)}$ is a diffeomorphism to its image. For any
$z\in Z$ denote $Z_z:=B(z,1) \cap Z$. For any $z \in Z$ there
exists $\eps>0$ such that $\nu|_{B_M(Z_z,\eps)}$ is a
diffeomorphism to its image, for example $\eps =
\min_{Z_z}(\sigma)$. Denote by $G(z)$ the supremum of all such
$\eps$. Let $\rho_Z^M(z) := min (G(z)/2, 1/10)$, $Tube(Z,M):=
\nu(B(Z,\rho_Z^M))$ and $\nu_Z^M := \nu|_{B(Z,\rho_Z^M)}$. It is
easy to see that $(Tube(Z,M), \rho_Z^M , \nu_Z^M)$ form a Nash
tubular neighborhood.

Clearly we can always find $\rho_Z^{M\setminus Y}< \rho_Z^{M}$
such that the corresponding Nash tubular neighborhood in $M$ does
not intersect $Y$. \proofend
\begin{corollary} \label{NashTubeCor}
Let $M$ be an affine Nash manifold and $Z \subset M$ be closed
affine Nash submanifold. Then $Z$ has a neighborhood $U$ in $M$
which is Nash diffeomorphic to the total space of the normal
bundle $N_Z^M$.
\end{corollary}
This corollary follows from the theorem using the following
technical lemma.
\begin{lemma}
Let $E \rightarrow M$ be a Nash bundle. Let $f\in {\cal N}(E)$ be
a Nash fiberwise homogeneous function of even degree $l \geq 2$.
Suppose $f(M)=\{0\}$ and $f(E \setminus M)>0$. Let $U:=\{x\in E|
f(x)<1\}$. Then $U$ is Nash diffeomorphic to $E$.
\end{lemma}
\emph{Proof.} Define a Nash diffeomorphism by stereographic
projection, i.e. $\nu:U \rightarrow E$ by
$\nu(m,v)=(m,\frac{1}{1-f(m,v)}v)$. \proofend
\section{Schwartz and tempered functions on affine Nash manifolds}
\label{AffSchwartz}

\subsection{Schwartz functions}
%
%\begin{definition}
%{Let $M$ be a smooth affine algebraic variety. We define
%\emph{the space of Schwartz functions on $M$} by $\Sc(M):=\{ f \in
%C ^ \infty (M) |$ for any algebraic differential operator $D$ on
%$M$ , $Df$ is bounded $\}$. We introduce topology on this space
%defined by the following system of semi-norms : $||f||_D := \sup
%\limits _ {x\in M} |Df(x)|$. It is easy to see that this topology
%is Frechet. Note that $\Sc(\R^n)$ is the space of classical
%Schwartz functions.}
%\end{definition}
%
\begin{definition}\label{das}
{Let $M$ be an affine Nash manifold. We define \textbf{the space
of Schwartz functions on $M$} by $\Sc(M):=\{ \phi \in C ^ \infty
(M) |$ for any Nash differential operator $D$ on $M$, $D\phi$ is
bounded $\}$. We introduce a topology on this space by the
following system of semi-norms : $||\phi||_D := \sup \limits _
{x\in M} |D\phi(x)|$. }
\end{definition}

%\newpage
\begin{proposition}\label{AlgCrit}
%In the case when $M$ is smooth affine algebraic variety, the 2
%definitions are equivalent. Moreover,
Let $M \subset \R^n$ be a closed affine Nash submanifold, and
$\phi$ be a smooth function on $M$. Suppose that for any algebraic
differential operator $D$ on $M$, $D\phi$ is bounded. Then
$\phi\in \Sc(M)$.
\end{proposition}

\emph{Proof.} Let $D'$ be a Nash differential operator. By lemma
\ref{multalg}, $D'\phi=\sum \limits _{i=1}^n g_iD_i\phi$ where
$D_i$ are algebraic differential operators on $M$ and $g_i$ are
Nash functions on $M$. By lemma \ref {majalg}, $g_i$ can be
majorated by polynomials $h_i$, hence $|D'\phi|=|\sum g_iD_i\phi|
\leq \sum \limits _{i=1}^n |(h_iD_i)\phi|$. The functions
$(h_iD_i)\phi$ are bounded since $h_iD_i$ is also an algebraic
differential operator. \proofend
\begin{corollary}
$\Sc(M)$ is a \Fre space. \footnote{This follows from the proof of
the previous proposition rather than from its formulation.}
\end{corollary}
\begin{corollary} \label{Usual}
If $M$ is a smooth algebraic variety then $\Sc(M)$ can be defined
as the space of smooth functions $\phi$ such that $D\phi$ is
bounded for every algebraic differential operator $D$ on $M$. In
particular, $\Sc(\R^n)$ is the space of classical Schwartz
functions.
\end{corollary}
\begin{remark}
It is easy to see that:\\
1. Nash differential operators act on Schwartz functions. \\
2. Schwartz functions form an algebra.
\end{remark}
\begin{remark}
One could give an equivalent definition of Schwartz functions
using a system of $L^2$ norms. Namely, choosing a Nash
non-vanishing density $\mu$ and defining $||\phi||_{D,\mu}:=
\int_M |D\phi|^2d\mu$.
\end{remark}
\subsection{Tempered functions}
\begin{definition} \label{Tempered}
{Let $M$ be an affine Nash manifold. A function $\alpha \in
C^{\infty}(M)$ is said to be \textbf{tempered} if for any Nash
differential operator $D$ there exists a Nash function $f$ such
that $|D\alpha| \leq f $.}
\end{definition}
\begin{proposition} \label{ModGrowth}
Let $M$ be an affine Nash manifold and $\alpha$ be a tempered
function on $M$. Then $\alpha \Sc(M)\subset \Sc(M)$.
\end{proposition}
\emph{Proof.} Let $\phi$ be a Schwartz function and $D$ be a Nash
differential operator on $M$. By the Leibnitz rule, $D(\alpha
\phi)= \sum \limits _{i=1}^n (D_i\alpha)(D'_i\phi)$ for some Nash
differential operators $D_i$ and $D'_i$. Since $\alpha$ is
tempered, there exist positive Nash functions $f_i$ such that
$|D_i\alpha| \leq f_i$. Denote $D''_i:= f_iD'_i$. So $|D(\alpha
\phi)|= |\sum \limits _{i=1}^n (D_i\alpha)(D'_i\phi)| \leq \sum
\limits _{i=1}^n |f_iD'_i\phi|=\sum \limits _{i=1}^n |D''_i\phi|$,
which is bounded since $\phi$ is Schwartz. \proofend
\begin{remark}
One can prove that the converse statement is also true, namely if
$\alpha \Sc(M)\subset \Sc(M)$ then $\alpha$ is tempered. We will
neither use nor prove that in this paper.
\end{remark}
The proof of the following lemma is straightforward.
%
%\begin{lemma}
%If $f$ is tempered then for any Nash differential operator $D$ and
%Schwartz function $g$, $gDf$ is bounded.
%\end{lemma}
%
\begin{lemma} \label{TempAlg}
1. Nash differential operators act on tempered functions.\\
2. Tempered functions form a unitary algebra and every tempered
function whose absolute value is bounded from below by a strictly
positive constant is invertible in this algebra.
\end{lemma}
\subsection{Extension by zero of Schwartz functions}
\begin{proposition} {(Extension by zero).} \label{ext0}
Let $M$ be an affine Nash manifold and $U \subset M$ be an open
(semi-algebraic) subset. Then the extension by zero of a Schwartz
function on $U$ is a Schwartz function on $M$ which vanishes with
all its derivatives outside $U$.
\end{proposition}
The proposition follows by induction from the following lemma.
\begin{lemma}
For any $\phi:U \rightarrow \R$ denote by $\widetilde{\phi}:M
\rightarrow \R$ its extension by 0 outside $U$. Let $\phi \in
\Sc(U)$. Then $\widetilde{\phi}$ is differentiable at least once
and for any Nash differential operator $D$ of order 1 on $M$,
$D\widetilde{\phi} = \widetilde{D|_U\phi}$.
\end{lemma}
\emph{Proof.} We have to show that for any $z\in M\setminus U$,
$\widetilde{\phi}$ is differentiable at least once at $z$ and its
derivative at $z$ in any direction is 0. Embed $M \hookrightarrow
\R^n$ and denote $F_z(x) := ||x-z||$. Clearly, $1/F_z^2 \in {\cal
N}(U)$. Hence $\phi/F_z^2$ is bounded in $U$ and therefore
$\widetilde{\phi} / F_z^2$ is bounded on $M \setminus z$, which
finishes the proof. \proofend
\subsection{Partition of unity} \label{AffPartUni}
In the proofs in this subsection we use the technical cover tools
developed in appendix (subsection \ref{Covers}).
\begin{theorem}{(Partition of unity).} \label{partuni}
Let $M$ be an affine  Nash manifold, and let $(U_i)_{i=1}^n$ be a
finite
open (semi-algebraic) cover of $M$. Then \\
1) there exist tempered functions $\alpha_1,\dots,\alpha_n$ on $M$
such that $supp(\alpha_i) \subset U_i$, $\sum \limits _{i=1}^n
\alpha_i
=1$. \\
2) Moreover, we can choose $\alpha_i$ in such a way that for any
$\phi \in \Sc(M)$, $\alpha_i \phi \in \Sc(U_i)$.
\end{theorem}
\emph{Proof of 1.} By lemma \ref{FinCor}, which follows from the
finiteness theorem, we can assume $U_i=\{x \in M| F_i(x)\neq 0 \}$
and $F_i|_{U_i}$ is a positive Nash function. By proposition
\ref{AffineNarCov}, there exists a strictly positive Nash function
$f$ such that the sets $V_i=\{x \in M| F_i(x)> f(x)\}$ also cover
$M$. Let $\rho: \R \rightarrow [0,1]$ be a smooth function such
that $\rho ((-\infty,0.1]) = \{0\}, \rho([1, \infty)) = \{1\} $.
Define $\beta_{i}:= \rho(F_{i}/f)$. It is easy to see that they
are tempered on $M$. Denote $\alpha_i := \frac{\beta_i}{\sum
\beta_i}$.The functions $\beta_i$ are tempered, $\sum \beta_i \geq
1$ hence $\alpha_i$ are tempered functions by lemma \ref{TempAlg}.
Hence the $\alpha_i$ form a tempered partition of unity. \proofend \\
In the proof of part 2 we will need the following technical lemma.
\begin{lemma} \label{techlemm}
Let $M\subset \R^d$ be an affine closed Nash submanifold and $$U=
\{x \in M|\, \forall 1 \leq i \leq n,\, p_{i}(x)>0 \}$$ where
$p_i$ are polynomials on $\R^d$. Let $g > 0$ be a Nash function
and $$U' = \{x \in M|\, \forall 1 \leq i \leq n, \, p_{i}(x)>g(x)
\}.$$ Then any Schwartz function $\phi$ on $M$ which is 0 outside
$U'$ is a Schwartz function on $U$.
\end{lemma}
\emph{Proof.} There is a standard closed embedding $\nu':U
\hookrightarrow \R^{d+n}$ whose last $n$ coordinates are defined
by $1/p_i$. By proposition \ref{AlgCrit}, it is enough to check
that for any differential operator $D$ on $U$ algebraic with
respect to $\nu'$, $D\phi|_U$ is bounded. Standard algebraic
geometry arguments show that $D$ is a sum of differential
operators of the form $\frac{1}{\prod p_i^k} (D'|_U)$ where $k$ is
a natural number and $D'$ is a differential operator on $M$
algebraic with respect to the imbedding $M \subset \R^d$. The
lemma follows now from the fact that inside $U'$, $\frac{1}{\prod
p_i} \leq g^n$. \proofend
\begin{remark}
In fact, more is true: any Schwartz function on $M$ which vanishes
outside $U$ together with all its derivatives is Schwartz on $U$.
We will prove that in subsection \ref{OpenChar}, but the proof
uses partition of unity.
\end{remark}
\emph{Proof of part 2 of theorem \ref{partuni}.} Fix a closed
imbedding $M \subset \R^d$. By finiteness theorem (\ref{finbas}),
we can suppose $U_i=\{x \in M| p_i(x) > 0 \}$ where $p_i$ are
polynomials on $\R^d$. By proposition \ref{AffineNarCov}, there
exists a Nash function $g>0$ such that the sets $V_i=\{x \in M|
p_i(x) > g \}$ also cover $M$. Let $\alpha_i$ be the tempered
partition of unity for the cover $V_i$. By the previous lemma,
$\alpha_i$ satisfy 2.
 \proofend
\begin{definition}
Let $M$ be an affine Nash manifold. We define \emph{the cosheaf
$\Sc_M$ of Schwartz functions on $M$} in the following way. For
any open (semi-algebraic) subset $U$ define $\Sc_M(U)$ to be
$\Sc(U)$ and for $V \subset U$ define the extension map
$ext_{V,U}:\Sc_M(V) \rightarrow \Sc_M(U)$ to be the extension by
zero.
\end{definition}

\begin{proposition} \label{affcosheaf}
$\Sc_M$ is a cosheaf (in the restricted topology).
\end{proposition}

\emph{Proof.} It follows from partition of unity (theorem
\ref{partuni}) and extension by zero (proposition \ref{ext0}).
\proofend%
\subsection{Restriction and sheaf property of tempered functions}
In the proofs given in this subsection we use the technical tools
developed in the Appendix (subsection \ref{SemiNash}).
\begin{theorem}{(Tempered restriction).} \label{TemRes}
Let $M\subset \R^n$ be an affine Nash manifold, $p_1, \dots, p_n$
be polynomials on $\R^n$ and $U=U_{p_1, \dots , p_n} = \{x \in M|
\, p_1(x)>0, \dots , p_n>0 \}$ be a basic open semi-algebraic
subset. Then the restriction to $U$ of any tempered function
$\alpha$ on $M$ is tempered.
\end{theorem}
\emph{Proof.} It is enough to prove the statement for $n=1$.
Denote $p=p_1$. Let $D$ be a Nash differential operator on $U$. We
will prove that
$|D\alpha|$ is bounded by some Nash function.\\
Step 1. Reduction to the case $D=f D'|_U$, where $D'$ is a Nash
differential operator on $M$ and $f$ is a Nash function on
$U$.\\
By lemma \ref{multalg}, we may assume that there exist an
algebraic differential operator $D_1$ and a Nash function $f_1$ on
$U$ such that $D=f_1 D_1$. From algebraic geometry we know that
there exist a natural number $n$ and an algebraic differential
operator $D'$ on
$M$ such that $D_1 = p^{-n}D'|_U$. Now we take $f=f_1 p^{-n}$.\\
Step 2. Proof of the theorem.\\
Since $\alpha$ is a tempered function on $M$, there exists a
positive Nash function $f_2$ on $M$ such that $|D'\alpha|\leq
f_2$. So $|D\alpha| \leq |f| \cdot f_2|_U $. By lemma
\ref{MajNash} there exists a Nash function $f_3$ on $U$ such that
$|f| \cdot f_2|_U \leq f_3$ and hence also $|D\alpha|\leq f_3$.
\proofend

\begin{theorem}{(The sheaf property of tempered functions).} \label{LocDef}
Let $M$ be an affine Nash manifold and $(U_i)_{i=1}^n$ its open
cover. Let $\alpha$ be a smooth function on $M$. Suppose
$\alpha|_{U_i}$ is tempered on $U_i$. Then $\alpha$ is tempered on
$M$.
\end{theorem}
\emph{Proof.} Embed $M \hookrightarrow \R^n$ and let $d$ be the
global metric induced from $\R^n$. Let $D$ be a Nash differential
operator on $M$. We know that there exist strictly positive Nash
functions $f_i$ on $U_i$ such that $|D\alpha| \leq f_i$ on $U_i$.
Let $G_i(x):=\min(f_i(x)^{-1}, d(x,M\setminus U_i))$. Extend $G_i$
by zero to $M$ and define $G:=\max G_i$. Let $F:=G^{-1}$ and let
$f$ be a Nash function that majorates $F$, which exists by lemma
\ref{MajNash}. It is easy to see that $|D\alpha|\leq f$. \proofend

\begin{definition}
Let $M$ be an affine Nash manifold. We define \emph{the sheaf
$\T_M$ of tempered functions on $M$} in the following way. For any
open (semi-algebraic) subset $U$ define $\T_M(U)$ to be $\T(U)$
and for $V \subset U$ define the usual restriction map
$res_{U,V}:\T_M(V) \rightarrow \T_M(U)$.
\end{definition}
We can summarize the previous 2 theorems in the following
proposition:
\begin{proposition} \label{affTemcosheaf}
$\T_M$ is a sheaf (in the restricted topology).
\end{proposition}
%
%This proposition follows from the last two theorems.
%
%\begin{construction} \label{SchwartzConst}
%Let $M$ be an affine Nash manifold and $Z \subset M$ be a closed
%affine Nash submanifold. We construct a function $f_Z^M:M
%\rightarrow \R $ that will be 1 in a neighborhood of $Z$ and 0
%outside another neighborhood of $Z$.
%
%Let $f: \R^{\geq 0} \rightarrow [0,1]$ be a smooth function such that
%$f([0,0.1])=\{1\}$ and $f([0.9,\infty)) = \{0\}$. Define $g:
%Tube(Z,M) \rightarrow \R $ by $g(z,v) = |v|/\rho_Z^M(z)$ and
%$fn_Z^M : N_{Z}^{M} \rightarrow \R $ by $fn_Z^M(z,v) = f(g(z,v))$
%if $(z,v) \in Tube_N(Z,M)$ and 0 otherwise. We define $f_Z^M : M
%\rightarrow \R$ by $f_Z^M = fn_Z^M \circ (\phi_Z^M)^{-1}$ on
%$Tube_N(Z,M)$ and $f = 0$ outside $Tube_N(Z,M)$.
%\end{construction}
\subsection{Extension of Schwartz and tempered functions from closed submanifolds}
\begin{theorem} \label{ext}
Let $M$ be an affine Nash manifold and $Z \hookrightarrow M$ be a
closed Nash submanifold. The restriction $\Sc(M) \rightarrow
\Sc(Z)$ is defined, continuous and \textbf{onto}. Moreover, it has
a section $s: \Sc(Z) \rightarrow \Sc(M)$ such that if $\phi\in
\Sc(Z)$ is 0 at some point $p$ with all its derivatives, then
$s(\phi)$ is also 0 at $p$ with all its derivatives.
\end{theorem}
%we will first prove a technical lemma.
%\begin{lemma}
%Let $M$ be a Nash manifold, $f \in \Sc(M)$, $\nu \in
%C^{\infty}(R^{\geq 0})$ satisfying $\nu([0.9,\infty) = 0$. Let $g$
%be a Nash function on $M$, $U := \{(x,v) \in M \times \R^n| |v| <
%g(x)\} $, $h$ be a Nash function on $U$, $F: U \rightarrow \R$
%defined by $F(x,v) =
%f(x)\nu(|v|/g(x))h(x,v)$.\\
%Then $F$ is Schwartz on $U$.
%
%Moreover, if all the derivatives of $\nu$ vanish at 0 (except maybe
%$\nu$ itself) and $f$ vanishes with all its derivatives at some
%point $x \in M$ then $F$ also vanishes with all its derivatives at
%$x$.
%
%\end{lemma}
%
%\emph{Proof}
%
%Step 1. $F$ is bounded.
%
%\emph{Proof}
% Denote $M=sup(\nu)$.  $\forall x \in M$, define $k(x) =
% \max_{|v| \leq 0.9}(h(x,v))$. By Tarski-Seidenberg theorem, it
%is semi-algebraic on $M$ and hence $fk$ is bounded. Denote its
%supremum by $N$. Hence $F \leq MN$.
%
%Step 2. $F$ is Schwartz.
%
%\emph{Proof} Every Nash differential operator $D$ on $U$ can be
%obtained from partial derivatives in $\R^n$, Nash differential
%operators on $M$ and multiplications by Nash functions on $U$ by
%taking compositions and sums. Each one of these operators
%preserves the conditions of the lemma, hence by step 1 $DF$ will
%be bounded, so $F$ is Schwartz. The last statement of the lemma
%follows from Leibnitz rule . \proofend \\
%%
\emph{Proof.} Clearly, the restriction is well defined and
continuous. Let us show that it is onto. \setcounter{case}{0}
\begin{case} $M=N\times \R^n$, $Z=N\times \{0\}$.\\
{\rm Choose a Schwartz function $\psi \in \Sc(\R^n)$ such that
$\psi=1$ in a neighborhood of the origin. For any $\phi\in \Sc(Z)$
we define $s(\phi)(n,v):=\phi(n,0)\psi(v)$. Clearly, $s$ is the
required section.}
\end{case}
\begin{case}
{\rm There exists a Nash diffeomorphism between $M$ and a Nash
vector
bundle over $Z$ which maps $Z$ to the zero section. \\
In this the claim follows from case 1 and theorem \ref{partuni}.}
\end{case}
\begin{case}
{\rm General.\\ Follows from case 2, corollary \ref{NashTubeCor}
of the Nash tubular neighborhood theorem and extension by zero
\ref{ext0}. \proofend}
\end{case}
\setcounter{case}{0}
In the same way one can prove an analogous theorem for tempered
functions.
\begin{theorem} \label{ExtTemp}
Let $Z \hookrightarrow M$ be a closed Nash embedding of affine
Nash manifolds. The restriction $\T(M) \rightarrow \T(Z)$ is
defined, continuous and onto. Moreover, it has a section $s: \T(Z)
\rightarrow \T(M)$ such that if $\alpha \in \T(Z)$ is 0 at some
point $p$ with all its derivatives, then $s(\alpha)$ is also 0 at
$p$ with all its derivatives.
\end{theorem}

\section{Schwartz, tempered and generalized Schwartz sections of Nash bundles over arbitrary Nash
manifolds} \label{SchAbsNash}
\subsection{Main definitions}
To define tempered and Schwartz functions on abstract Nash
manifolds we use these notions on affine Nash manifolds and glue
them using sheaf and cosheaf properties respectively.
\begin{definition}
{Let $M$ be a Nash manifold, and $E$ be a Nash bundle over it. Let
$M=\bigcup \limits _{i=1}^k U_i$ be an affine Nash trivialization
of $E$. A global section $s$ of $E$ over $M$ is called
\textbf{tempered} if for any $i$, all the coordinate components of
$s|_{U_i}$ are tempered functions. The space of global tempered
sections of $E$ is denoted by $\T(M,E)$. }
\end{definition}
%\begin{definition}
%{let $M$ be a Nash manifold, and $E$ be a Nash bundle over it.
%A global section $s$ of $E$ over $M$ is called \emph{tempered} if
%for any open affine Nash submanifold $U\subset M$ and for any
%trivialization of $E|_U$, $s|_U$ is a tuple of tempered functions.
%The space of global tempered sections of $E$ is denoted by
%$\T(M,E)$. }
%\end{definition}
%
\begin{proposition}
The definition does not depend on the choice of $U_i$.
\end{proposition}
\emph{Proof.} It follows from the sheaf properties of tempered
functions (proposition \ref{TemRes} and theorem \ref{LocDef}).
\proofend
\begin{definition}\label{dlsglob}
Let $M$ be a Nash manifold, and let $E$ be a Nash bundle over it.
Let $M=\bigcup \limits _{i=1}^k U_i$ be an affine Nash
trivialization of $E$. Then we have a map $\phi : \bigoplus
\limits _{i=1}^k \Sc(U_i)^n \rightarrow C^\infty (M,E) $. We
define \textbf{the space $\Sc(M,E)$ of global Schwartz sections of
$E$} by $\Sc(M,E):= Im\phi $. We define the topology on this space
using the isomorphism $\Sc(M,E) \cong \bigoplus \limits _{i=1}^k
\Sc(U_i)^n / Ker\phi$.
\end{definition}
\begin{proposition}
The definition does not depend on the choice of $U_i$.
\end{proposition}
\emph{Proof.} It follows from theorems on partition of unity
(theorem \ref{partuni}) and extension by zero (proposition
\ref{ext0}). \proofend
\begin{definition}\label{dls}
{Let $M$ be a Nash manifold, and let $E$ be a Nash bundle over it.
We define \textbf{the cosheaf $\Sc_M^E$ of Schwartz sections of
$E$ } by $\Sc_M^E(U):= \Sc(U,E|_U)$. We also define \textbf{the
sheaf $\G_M^E$ of generalized Schwartz sections of $E$ } by
$\G_M^E(U):= (\Sc_M^{\widetilde{E}})^*$, where $\widetilde
{E}=E^*\otimes D_M$, and \textbf{the sheaf $\T_M^E$ of tempered
sections of $E$} by $\T_M^E(U):= \T(U,E|_U)$.  }
\end{definition}
\begin{proposition}
$\Sc_M^E$ is a cosheaf, and $\T_M^E$ and $\G_M^E$ are sheaves.
\end{proposition}
\emph{Proof.} It follows from the definitions. \proofend
\begin{remark}
$\T_M$  and $\G_M$ are functors from the category of Nash bundles
over $M$ to the categories of sheaves on $M$. Also, $\Sc_M$, is a
functor from the category of Nash bundles over $M$ to the
categories of cosheaves on $M$. The mappings of morphisms are
obvious.
\end{remark}
\subsection{Partition of unity} \label{SecPartUni}
This subsection is rather similar to subsection \ref{AffPartUni}
about partition of unity for affine Nash manifolds. In particular,
in the proofs in this subsection we use the technical cover tools
and definitions given in the Appendix (subsection \ref{Covers}).
\begin{theorem}{(Partition of unity for any Nash manifold).}
\label{anypartuni} Let $M$ be a Nash manifold, and let
$(U_i)_{i=1}^n$ be a finite open cover. Then \\
1) there exist tempered functions $\alpha_1,...,\alpha_n$ on $M$
such that $supp(\alpha_i) \subset U_i$, $\sum \limits _{i=1}^n
\alpha_i
= 1$. \\
2) Moreover, we can choose $\alpha_i$ in such a way that for any
$\phi \in \Sc(M)$, $\alpha_i \phi \in \Sc(U_i)$.
\end{theorem}

\emph{Proof of 1.}

Let $\{F_j\}$ be a basic collection of continuous semi-algebraic
functions such that $M_{F_j}$ is a refinement of $\{U_i\}$ . It
exists by proposition \ref{FinCov}. Let $\rho: \R \rightarrow
[0,1]$ be a smooth function such that $\rho ((-\infty,0.1]) =
\{0\}, \rho([1, \infty)) = \{1\} $. Denote $\beta _j:=\rho \circ
F_j$ and $\gamma_j = \frac{\beta_j}{\sum \beta_j}$. It is easy to
see that $\gamma_j$ are tempered. Now for every $j$ we choose
$i(j)$ such that $M_{F_j} \subset U_{i(j)}$. Define
$\alpha_i:=\sum \limits _{j|i(j)=i}\gamma_j$. It is easy to see
that $\alpha_j$ is a tempered partition of unity.
%$V_j,V'_j,V''_j$ be finite open (semi-algebraic) covers of $M$ and
%$f_j:V_j \rightarrow \R^{\geq 0}$ be continuous semi-algebraic
%functions such that $V_j$ are affine, the cover $V_j$ is a
%refinement of $U_i$, the cover $V'_j$ is a proper refinement of
%the cover $V_j$, $V'_j=\{x\in V_j|f_j(x) \neq 0\}$, $f_j|_{V_j'}$
%is Nash positive function and $V''_j=\{x \in V_j|f_j(x)>1\}$. They
%exist by proposition \ref{FinCov} from the appendix.
%
%It is enough to prove the theorem for the cover $V_j$. Let $\rho:
%\R \rightarrow [0,1]$ be a smooth function such that $\rho
%((-\infty,0.1]) = \{0\}, \rho([1, \infty)) = \{1\} $. Denote
%$\beta _j:=\rho \circ f_j$ continued by 0 outside $V_j$. It is
%easy to see that $\alpha_j = \beta_j/(\sum \beta_j)$ gives a
%tempered partition of unity.
\proofend\\
As in the proof of affine partition of unity theorem, part 2
follows from part 1, proposition \ref{NarCov} and the following
lemma, similar to lemma \ref{techlemm}.
\begin{lemma} \label{Blin}
Let $M$ be a Nash manifold. Let $V\subset U \subset M$ be open
(semi-algebraic) sets such that $\overline{V} \subset U$. Let
$\alpha$ be a tempered function on $M$ supported in $V$. Let
$\phi$ be a Schwartz function on $M$. Then $\alpha \phi|_U$ is a
Schwartz function on $U$.
\end{lemma}
In order to prove this lemma we need the following technical
lemma.
\begin{lemma}
Let $M$ be an affine Nash manifold. Let $\phi$ be a Schwartz
function on $M$. Let $U \subset M$ be an open (semi-algebraic)
subset. Denote $\phi':=\phi \mathbf{1}_U$ where $\mathbf{1}_U$ is
the characteristic function of $U$. Suppose that $\phi'$ is
smooth. Then $\phi'$ is Schwartz.
\end{lemma}
\emph{Proof.} Denote $W :=U \cup (\overline{U})^c$. Let $D$ be a
Nash differential operator on $M$. It is easy to see that
$|D\phi'|_W|\leq |D\phi|_W|$ and $W$ is dense in $M$. Hence
$|D\phi'| \leq |D\phi|$ which is bounded. \proofend

\emph{Proof of lemma \ref{Blin}.} Without loss of generality we
can assume that $M$ is affine. Denote $W:=M \setminus
\overline{V}$. Note that $M=U \cup W$. Since Schwartz functions
form a cosheaf, $\phi=\phi_U + \phi_W$ where $\phi_U$ and $\phi_W$
are extensions by zero of Schwartz functions on $U$ and on $W$,
respectively. Note that
$\phi|_{\overline{V}}=\phi_U|_{\overline{V}}$ so $\phi|_U =
(\phi_U)|_U \cdot(\mathbf{1}_{\overline{V}})|_U$ and therefore by
the lemma
$\phi|_U$ is Schwartz. \proofend%
%
%Let $U_i$ be an affine cover of $M$. It is enough to prove for $g$
%which is an extension by zero of a Schwartz function on $U_i$.
%Consider the following cover of $U_i$: $U_i= (U_i \setminus
%\overline{V}) \cup (U_i \cap U)$. By partition of unity for affine
%case (theorem \ref{partuni}), $g$ can be represented as a sum $g =
%g_1 + g_2$ where $g_1$ is extension by zero of a Schwartz function
%on $U_i \setminus \overline{V}$ and $g_2$ is extension by zero of
%a Schwartz function on $(U_i \cap U)$. In particular, $g_2$ is an
%extension by zero of a Schwartz function on $U_i$. So $fg = fg_2$
%and its restriction to $U$ is Schwartz.
%\proofend \\
%
%\emph{Proof of partition of unity.} Suppose without loss of
%generality that $U_i$ are affine. Choose an open cover
%$(V_i)_{i=1}^n$ such that $V_i \subset U_i$ for all i. It is possible
%by proposition \ref{narcov}. By finiteness theorem (\ref{finbas})
%$V_i = \bigcup V_{ij}$ where $V_{ij}=\{x\in U_i|\quad
%f_{ijk}(x)>0\}$. From this point the proof is very similar to the
%proof for affine case (theorem \ref{partuni}) but using the
%previous lemma instead of lemma \ref{techlemm}. \proofend \\

\begin{remark}
Partitions of unity clearly exist for Schwartz sections,
generalized Schwartz sections and tempered sections of any Nash
bundle, i.e. for any Nash bundle $E$ over $M$ and any finite open
cover $U_i$ of $M$ there exists a tempered partition of unity
$\alpha_i$ such that for any $\phi \in \Sc(M,E)$, $\phi\alpha_i
\in \Sc(U_i,E)$ and for any $\xi \in \G(M,E)$, $\xi\alpha_i \in
\G(U_i,E)$.
\end{remark}
\subsection{Basic properties} \label{BasProp}
Let us now prove properties \ref{p1} - \ref{p5} mentioned in
section \ref{Summary}. \\
Property \ref{p1} holds by definition.
Clearly, for affine $M$ and trivial 1-dimensional $E$, $\Sc(M,E) =
\Sc(M)$ and property \ref {p2} holds. Property \ref {p4} is
satisfied by definition.
%Property \ref {p4'} easily follows from proposition \ref{ModGrowth}.
Property \ref {p5} follows from theorems \ref{ext} and
\ref{ExtTemp}.
\begin{theorem} \label{ProofProp3}
Property \ref{p3} holds, i.e. for compact Nash manifold $M$,
$\Sc(M,E) = \T(M,E)=C^{\infty}(M,E)$.
\end{theorem}
 \emph{Proof} Let $\alpha$ be a
smooth section of a Nash bundle $E$ over a compact Nash manifold
M. We have to show that it is also a Schwartz section. Since $M$
is a Nash manifold, any point $m \in M$ has a neighborhood $U_m$
Nash diffeomorphic to an open ball in $\R^n$ of radius 1. Denote
by $V_m \subset M$ the preimage of the ball of radius 0.9 under
this Nash diffeomorphism. $\{{V_m}\}_{m \in M}$ is a cover of $M$.
Let us choose a finite subcover $V_{m_i}$. By classical partition
of unity, $\alpha$ can be represented as $\alpha=\sum \alpha_i$
where $supp(\alpha_i)\subset V_{m_i}$. $\alpha_i$ are clearly
Schwartz on $U_i$. Hence $\alpha$ is Schwartz on $M$. \proofend

\subsection{Characterization of Schwartz functions on open subset}
\label{OpenChar}

Let us now prove property \ref {p6}.
First we will prove it for trivial 1-dimensional bundle.

\begin{theorem}\label{opensetTriv} {(Characterization of Schwartz
functions on open subset)}\\ Let $M$ be a Nash manifold, $Z$ be a
closed (semi-algebraic) subset and $U= M\setminus Z$. Let $W_Z$ be
the closed subspace of $\Sc(M)$ defined by $W_Z:=\{\phi \in \Sc(M)
|\phi$ vanishes with all its derivatives on $Z\}$. Then restriction
and extension by 0 give an isomorphism $\Sc(U) \cong W_Z$.
\end{theorem}

We will use the following elementary lemma from calculus.
\begin{lemma}\label{ThR}
Suppose $\alpha \in C^{\infty}(\R)$ vanishes at 0 with all its
derivatives. Then for any natural number $n$, $\alpha(t) = t^n
\alpha^{(n)}(\theta)$ for some $\theta \in [0,t]$.
\end{lemma}

\emph{Proof of the theorem.}

\setcounter{case}{0}
\begin{case}
{\rm
$M=\R^N$.\\
Let $\phi \in \Sc(U)$ and let $\widetilde{\phi}$ be its extension
by 0. By proposition \ref{ext0} on extension by zero,
$\widetilde{\phi} \in W$.

Now, let $\phi \in W_Z$. For any point $x \in \R^N$ define
$r(x):=dist(x,Z)$. Let $S:=S(0,1) \in \R^N$ be the unit sphere.
Consider the function $\psi$ on $S \times Z \times \R$ defined by
$\psi(s,z,t):=\phi(z+ts)$. From the previous lemma \ref{ThR} we
see that $\psi(s,z,t)=t^n \frac{\partial^n }{(\partial t)^n} \psi
(x,s,t)|_{t=\theta}$ for some $\theta \in [0,t]$. As $\phi$ is
Schwartz, it is easy to see that $\frac{\partial^n }{(\partial
t)^n} \psi (x,s,t)$ is bounded on $Z \times S \times \R$.
Therefore $|\psi(s,z,t)|\leq C |t|^n$ for some constant $C$ and
hence $\phi/r^n$ is bounded on $\R^N$ for any $n$.

Let $h$ be a Nash function on $U$. By lemma \ref{saloc}, $r^nh$
extends by 0 to a continuous semi-algebraic function on $\R^N$ for
$n$ big enough. It can be majorated by $f \in {\cal N}(\R^N)$.
Therefore $|\phi h| = |(\phi/r^n)r^nh| \leq |\phi f|/r^n$. $\phi
f\in W$, thus $|\phi f|/r^n$ is bounded and hence $|\phi h|$ is
bounded.

For any Nash differential operator $D$ on $\R^N$, $D\phi \in W$.
Hence $hD\phi$ is bounded. By lemma \ref{DiffonOpen}, every Nash
differential operator on $U$ is a sum of differential operators of
the form $hD|_U$, where $D$ is a Nash differential operator on
$\R^N$ and $h$ a Nash function on $U$. Hence $\phi|_U \in
\Sc(U)$.}
\end{case}
%\item
%\begin{case} $M$ is affine Nash $M$, Z is compact.
%
%Follows from the previous step and theorem \ref{ext} (extension from
%closed Nash submanifold).
%\end{case}

%\item
%\begin{case} $M$ is semi-algebraic open subset of $S^n$ (and $Z$ is
%arbitrary) \\
%Follows trivially from the previous step.
%\end{case}

%\begin{case}%\item
%$M=\R^N$ \\
%Consider the stereographic embedding $\R^N \hookrightarrow S^N$.
%\end{case}

\begin{case}%\item
{\rm $M$ is affine. \\Follows from the previous case and theorem
\ref{ext} (extension from a closed Nash submanifold).}
\end{case}
%\item
\begin{case}
{\rm General case.\\
Choose an affine cover of $M$. The theorem now follows from the
previous case and partition of unity. \proofend}
\end{case}
\setcounter{case}{0}

Property \ref{p6} is an immediate corollary of the previous
theorem and partition of unity. Let us remind it.

\begin{theorem}\label{openset} {(Characterization of Schwartz
sections on open subset)}\\ Let $M$ be a Nash manifold, $Z$ be a
closed (semi-algebraic) subset and $U= M\setminus Z$. Let $W_Z$ be
the closed subspace of $\Sc(M,E)$ defined by $W_Z:=\{\phi \in
\Sc(M,E) |\phi$ vanishes with all its derivatives on $Z\}$. Then
restriction and extension by 0 give an isomorphism $\Sc_M^E(U)
\cong W_Z$.
\end{theorem}

%\begin{theorem}\label{openset} {(Characterization of Schwartz
%sections on open subset)}\\ Let $M$ be a Nash manifold, $E$ a Nash
%bundle over $M$, $U \subset M$  an open (semi-algebraic) subset ,
%and $Z=M\setminus U$. Then restriction and extension by 0 give an
%isomorphism $\Sc_M^E(U) \cong \{\phi \in \Sc(M,E) |\phi$ vanishes
%with all its derivatives on $Z\}$.
%\end{theorem}

\begin{corollary} \label{flabbiness}
Let $E \rightarrow M$ be a Nash bundle. Then the cosheaf $\Sc_M^E$
and the sheaf $\G_M^E$ are flabby. In other words, let $V \subset
U \subset M$ be open (semi-algebraic) subsets. Then
$ext_{V,U}:\Sc_M^E(V) \rightarrow \Sc_M^E(U)$ is a closed
embedding and $res_{U,V}:\G_M^E(U) \rightarrow \G_M^E(V)$ is onto.
\end{corollary}

\emph{Proof.} The map $ext_{V,U}$ is a closed embedding by the
theorem. Hence by the Hahn Banach theorem the dual map $res_{U,V}$
is onto. \proofend

\subsubsection{Remarks}
\begin{remark}
Some of the ideas of the proof of theorem \ref{opensetTriv} are
taken from Casselman's unpublished work \cite{Cas2}.
\end{remark}

%\begin{remark}
%The proof of case 1 can give a proof of case 3 if it is done more
%delicately. It would make cases 2 and 3 obsolete.
%\end{remark}

\begin{remark}
Since extension by zero $\Sc(U) \cong W_Z$ is a continuous linear
isomorphism, the inverse map is also continuous by Banach open map
theorem. In fact, as it often happens, our proof that extension is
onto can be easily refined to prove that the inverse map is
continuous.
\end{remark}

\subsection{Generalized Schwartz sections supported on closed
submanifolds}

\begin{definition}\label{suppdist}
{Let $M$ be a restricted topological space, and $F$ be a sheaf on
$M$. Let $Z \subset M$ be a closed subset. A global section of $F$
is said to be \textbf{supported in $Z$} if its restriction to the
complement of $Z$ is zero.}
%Let $Z \subset M$ be a closed subset and $U$ be its complement. We
%define the space $\Gamma_Z(M,F)$ of global sections of $F$
%supported in $Z$ by $\Gamma_Z(M,F):= Ker(res_{M,U})$
\end{definition}
\begin{remark}
Unfortunately, if we try to define support of a section, it will
not be a closed set in general, since infinite intersection of
closed sets in the restricted topology does not have to be closed.
Also, we cannot in general consider its closure, because there is
no closure in restricted topology by the same reason.
\end{remark}
%{Let $\underset {M} { \overset {E}{\downarrow}}$ be a Nash
%bundle. Let $\xi \in G\Sc^E_M(M)$. we will say that \emph{$\xi$ is
%supported in $Z$} iff for any $f \in \Sc^{\widetilde{E}}_M(M)$ which
%vanishes with all its derivatives on $Z$, $\xi (f) =0$.}
%\begin{remark}
%
%Using technic similar to the sleeve technic used in theorem
%\ref{openext}, it can be shown that this definition is equivalent to
%the following one:
%
%Let $\underset {M} { \overset {E}{\downarrow}}$ be a Nash bundle.
%Let $\xi \in \G^E_M(M)$. we will say that $\xi$ is supported in
%$Z$ iff for any $f \in \Sc^{E^*\otimes D_M}_M(M)$ which vanishes in
%an open semi-algebraic neighborhood of $Z$, $\xi (f) =0$.
%\end{remark}
\begin{proposition}
Let $E \rightarrow M$ be a Nash bundle and $Z \subset M$ be a
closed (semi-algebraic) subset. Let $\xi \in \G^E_M(M)$. Then
$\xi$ is supported in $Z$ if and only if for any $\phi \in
\Sc^{\widetilde{E}}_M(M)$ which vanishes with all its derivatives
on $Z$, $\xi (\phi) =0$.
\end{proposition}

\emph{Proof.} It is an immediate corollary of the characterization
of Schwartz sections on an open subset (theorem
\ref{openset}).\proofend

\begin{lemma}\label{l1}
Let $M$ be a Nash manifold and let $Z\subset M$ be a closed Nash
submanifold. Let $\Delta$ be the map $$ \{\phi \in \Sc(M) |\phi
\text{ is 0 on } Z \text{ with first i-1 derivatives} \}
\rightarrow \Sc(Z,S^i(CN_Z^M)),$$ given by the $i$-th derivative,
where $S^i$ means i-th symmetric power and $CN_Z^M$ is the
conormal bundle. Then $\Delta$ is well defined and onto.
\end{lemma}

\emph{Proof.} For $M = Z \times \R^d$ the lemma is trivial. For
the general case, it is proved in the same way as theorem
\ref{ext} - using Nash tubular neighborhood. \proofend

\begin{corollary}\label{c4}                                             %??? Yesh Taut !!!
Property \ref {p7} holds. Namely, let $Z \subset M$ be a closed Nash       %It seams to be cotangent bundle
subset. Consider $V=\{\xi \in \G(M,E) |\xi$ is supported in $Z
\}$. It has a canonical filtration $V_i$ such that its factors are
canonically isomorphic to $\G(Z,E|_Z \otimes S^i(N_Z^M) \otimes
D_M^*|_Z \otimes D_Z)$ where $N_Z^M$ is the normal bundle.
\end{corollary}

\emph{Proof}. Denote $K:= \Sc_M^{\widetilde{E}}(M)$,
$L:=\Sc^{\widetilde{E}}_M(U)$. Then $V =(K/L)^*$. Define
$K_i:=\{\phi \in K|$ for any Nash differential operator $D$ of
degree $\leq i-1$, $D\phi|_Z = 0\}$. $K_0=K$ and $\bigcap K_i = L$
and by the previous lemma $K_i/K_{i+1}=\Sc(Z,\widetilde{E}|_Z
\otimes S^i(CN_Z^M)) $. Define $V_i = (K/K_i)^*$. It is the
requested filtration.  \proofend \\
This corollary appeared in a similar form in Casselman's
unpublished work \cite{Cas2} and in another similar form in
\cite{CHM}.
\appendix
\section{Appendix} \label{App}
\subsection{Nash structures on standard bundles} \label{SecBunApp}
In this section we construct Nash structures on standard bundles
from differential geometry.
\begin{theorem}
The tangent bundle of any Nash manifold has a canonical structure
of a Nash bundle.
\end{theorem}

\begin{corollary} \label{BunApp}
Tangent, normal and conormal bundles, the bundle of differential
$k$-forms, etc., of a Nash manifold have canonical structure of
Nash bundles.
\end{corollary}

\emph{Proof of the theorem.}
In order to prove this theorem we have to prove two statements:\\
1. The total space of the tangent bundle has a canonical
structure of a Nash manifold.\\
2. There is a finite cover of $M$ by Nash trivializations of the
tangent bundle.

%We are ready now to prove 1, the proof of 2 will be given later.
Proof of 1. It is enough to prove the proposition for a closed
affine Nash submanifold $M\subset \R^n$. In this case, the tangent
bundle of $M$ is defined in $\R^{2n}$ by semi-algebraic
conditions. Hence, by the Tarski-Seidenberg principle it is a
closed affine Nash submanifold of $\R^{2n}$. So the proposition
holds true for the tangent bundle, and hence also for normal and
conormal bundles, the bundle of differential $k$-forms, etc.

Proof of 2. It is sufficient to prove for the cotangent bundle
instead of the tangent one. Also, we can suppose that the manifold
is affine.\\
Let $M \subset \R^n$ be a closed Nash embedding. We have a map
$\pi$ from $(\R^n)^* $ to covector fields on $M$.
%By Nash covector fields we mean Nash maps from $M$ to the total space
%of the cotangent bundle that are sections of the natural
%projection.
Let $\{e_i\}$ be a basis of $(\R^n)^* $, $d:=\dim M$. For any
subset $S \subset \{1,\dots ,n\}$ of cardinality $d$ we define
$f_S:=\det(\pi(e_i), i \in S )$. Let $U_S:=\{x \in M | f_S(x) \neq
0\}$. Since $\pi$ is fiberwise onto, we know from linear algebra
that $U_S$ cover $M$. Clearly, every $U_S$ is an open
semi-algebraic subset and the cotangent bundle is trivializable on
it. \proofend

\begin{proposition} \label{BunAppCor}
Let $M\subset \R^n$ be a closed Nash submanifold. Then the space
of Nash covector fields on $M$ is generated over ${\cal N}(M)$ by
the covector fields $dx_i$.
\end{proposition}

\emph{Proof.} The fields $dx_i$ are exactly the fields $\pi(e_i)$
from the proof of the last theorem. \proofend

\subsubsection{The bundle of densities $D_M$} \label{BunDenApp}
Now we wish to define structure of Nash bundle on the bundle of
densities. It is done in the same way as in differential geometry.

We remind that the $\Z / 2\Z$- torsor $Orient_M$ of orientations
is defined as the quotient of the bundle $\Omega^{top}_M$ of top
differential forms by the multiplicative action of $\R_{>0}$. The
standard correspondence between torsors and bundles gives us a
bundle of orientations $BOr = Orient_M \bigotimes \limits _{\Z /
2\Z} \R$.

The bundle of densities $D_M$ is defined as the tensor product of
$BOr_M$ and $\Omega^{top}_M$, and its global smooth sections are
smooth measures on $M$.
%For more detailed description cf. ???.
%
\begin{definition}
{Let $M$ be a Nash manifold. We would like to define a Nash
structure on the bundle $BOr_M$ of orientations on $M$. Consider
the $\Z / 2\Z $ torsor $ Orient_M$ of orientations on $M$. Since
$\Omega^{top}_M$ is a Nash bundle, there exists a finite cover of
$M$ by open (semi-algebraic) subsets $M= \bigcup U_i$ such that
$Orient_M|_{U_i}$ is isomorphic to $\Z / 2\Z \times U_i$ as a
smooth locally trivial fibration. Choose such isomorphisms
$\beta_i$. These $\beta_i$ define the structure of a Nash locally
trivial fibration on $Orient_M$. It does not depend on the choice
of $\beta_i$ since the fiber is finite. We define \textbf{the
bundle $BOr_M$ of orientations on $M$} to be the bundle
corresponding to the $\Z / 2\Z $-torsor $ Orient_M$ and the sign
representation of $\Z / 2\Z $. In other words, $BOr_M = (Orient_M
\times \R)/(\Z / 2\Z)$ where $\Z / 2\Z$ acts diagonally, and on
$\R$ it acts by sign. It has an obvious structure of a Nash
bundle.}
\end{definition}
\begin{definition}
{We define \textbf{the bundle $D_M$ of densities on $M$} by $D_M =
BOr_M \bigotimes \limits _M \Omega^{top}_M$.}
\end{definition}
\subsection{Semi-algebraic notions on Nash manifolds}
\label{SemiNash}
\begin{definition}
{A subset $A$ of an affine Nash manifold $M\subset\R^n$ is called
\textbf{semi-algebraic} iff it is semi-algebraic in $\R^n$.
Clearly this notion does not depend on the embedding
$M\hookrightarrow \R^n$.

A subset $A$ of a Nash manifold $M$ is called
\textbf{semi-algebraic} iff its intersection with any open affine
Nash submanifold is semi-algebraic .

A map $\nu$ between Nash manifolds $M$ and $N$ is called
semi-algebraic iff its graph is a semi-algebraic subset of
$M\times N$.}
\end{definition}
We will need the following technical lemma, which is a direct
consequence of lemma \ref{majalg} about majoration of
semi-algebraic functions by polynomials.
\begin{lemma} \label{MajNash}
Let $M$ be an affine Nash manifold. Then any continuous
semi-algebraic function on it can be majorated by a Nash function,
and any continuous strictly positive semi-algebraic function on it
can be bounded from below by a strictly positive Nash function.
\end{lemma}
\begin{remark}
Let $\phi:M \to N$ be a semi-algebraic map between Nash manifolds.
Clearly, $\phi$ is continuous as a map of classical topological
spaces if and only if $\phi$ is continuous as a map of restricted
topological spaces.
\end{remark}
\begin{theorem}
Let $M$ be a Nash manifold. Then there exists a semi-algebraic
continuous embedding $\nu:M \hookrightarrow \R^n$, that means a
semi-algebraic map $\nu:M \rightarrow \R^n$ which is a
homeomorphism to its image.
\end{theorem}
For proof see \cite{Shi}, page 142, theorem III.1.1 .
\begin{corollary} \label{Met}
Let $M$ be a Nash manifold. Then there exists a semi-algebraic
continuous metric $d:M\times M \rightarrow \R$.
\end{corollary}
\begin{remark}
This corollary can be proven directly by defining the metric for
open affine Nash subsets that cover $M$ and then gluing.
\end{remark}
\begin{proposition}
Let $U$ be an open (semi-algebraic) subset of a Nash manifold $M$.
Then its closure $\overline{U}$ in the usual topology is closed
semi-algebraic, that is its complement $\overline{U}^c$ is open
semi-algebraic.
\end{proposition}
\emph{Proof.} It is enough to prove for affine $M$. This case
easily follows from Tarski-Seidenberg principle. \proofend
\begin{remark}
The last proposition shows that open sets have closure in the
restricted topology.
\end{remark}

\subsection{Covers} \label{Covers}
In this section we give some definitions and propositions that
help us to work with covers. The propositions of this subsection
are versions of a general statement which says that any open cover
of a Nash manifold can be replaced by a finer cover that has some
nice properties. This subsection is used in the proofs of
partition of unity both for affine and general Nash manifolds
(sections \ref{AffPartUni} and \ref{SecPartUni}). The central
statement of this section is proposition \ref{FinCov}.
\begin{notation}
Let $M$ be a Nash manifold and $F$ be a continuous semi-algebraic
function on $M$. We denote $M_F:=\{x \in M |F(x) \neq 0\}$.
\end{notation}
\begin{definition}
Let $M$ be a Nash manifold. A continuous semi-algebraic function
$F$ on $M$ is called \textbf{basic} if $F|_{M_F}$ is a positive
Nash function.
\end{definition}
The the above definition is motivated by the following lemma.
\begin{lemma} \label{FinCor}
Let $M$ be an affine Nash manifold. Then it has a basis of open
sets of the form $M_F$ where $F$ is a basic function.
\end{lemma}
This lemma follows directly from the finiteness theorem
(\ref{finbas}).
\begin{definition}
{We say that a cover $M=\bigcup \limits _{j=1}^m V_j$ is a
\textbf{refinement} of the cover $M=\bigcup \limits _{i=1}^n U_i$
if for any $j$ there exists $i$ such that $V_j \subset U_i$. \\
We say that a cover $M=\bigcup \limits _{j=1}^m V_j$ is a
\textbf{proper refinement} of the cover $M=\bigcup \limits
_{i=1}^n U_i$ if for any $j$ there exists $i$ such that
$\overline{V_j} \subset U_i$. }
\end{definition}
%
%\begin{definition}
%{\rm Let $U\subset V$ be open (semi-algebraic) subsets of a Nash
%manifold $M$. We say that $U$ is a proper }
%\end{definition}
%
\begin{proposition} \label{AffineNarCov}
Let $M=  \bigcup \limits _{i=1}^n U_i$ be a finite open
(semi-algebraic) cover of an affine Nash manifold $M$ such that
$U_i=\{x\in M| F_{ij}(x)>0$ for $j=1 \dots n_i\}$ for certain
continuous semi-algebraic functions $F_{ij}:M \rightarrow \R$.
Then there exists a strictly positive Nash function $g$ such that
the sets $V_i:=\{x\in M| F_{ij}(x)>g(x)$ for $j=1 \dots n_i\}$
also cover $M$.
\end{proposition}
\emph{Proof.}\\ Denote $G_{ij}:=max(F_{ij},0), \, G' :=
\frac{1}{2} \max \limits _{i=1}^n \min \limits _{j=1}^n G_{ij}$.
By lemma \ref{MajNash} there exists a Nash function $g$ such that
$0<g \leq G'$. \proofend
\begin{proposition} \label{NarCov}
Let $M=\bigcup \limits _{i=1}^n U_i$ be a finite open
(semi-algebraic) cover of a Nash manifold $M$. Then there exists a
finite open (semi-algebraic) cover $M=\bigcup \limits _{i=1}^n
V_i$ which is a proper refinement of $\{U_i\}$.
\end{proposition}
%In order to prove this proposition we will need the following
%lemma.
%\begin{lemma}
%Let $M$ be a Nash manifold. Then there exists a semi-algebraic
%metric $d:M\times M \rightarrow \R$ on $M$ such that the topology
%defined by $d$ coincides with the classical topology.
%\end{lemma}
%In fact, this lemma is not concerned with the smooth structure,
%and follows from the theorem that says that every Nash manifold
%can be $C^r$-Nash embedded into $\R^n$
%\emph{Proof.}\\
%Case 1 $M$ is affine.\\
%We embed $M\subset \R^n$ and take $d$ to be the Euclidian
%distance.\\
%Case 2 General. \\
%Choose an open affine cover $M=\bigcup \limits _{i=1}^n U_i$ and
%metrics $\mu_i$ on $U_i$. Define equivalence relation on $M$ by $x
%~ y$ iff there exists a finite sequence $x=x_0,x_1,\dots,x_k=y$ such that
%$x_i$ and $x_{i+1}$ lie in the same $U_{l_i}$. Note that if $x ~
%y$ then we can assume that there exists such a sequence of length
%$n$ and that the number of equivalence classes is at most $n$. Choose a
%system of representatives $r_1,\dots,r_l$. For every equivalence
%class we define $d_k(x,y) = \inf_{\{x=x_0,x_1,\dots,x_n=y|x_i \text{
%and }x_{i+1} \text{ lie in the same } U_{l_i} \}} \sum _{i=1}^n
%\mu_{l_{i-1}}(x_{i-1},x_{i})$. Now we define the metric $d$ by
%$d(x,y) = d_i(x,y)$ if $x ~ y ~ r_i$ and
%$d(x,y)=d_i(x,r_i)+1+d_j(r_j,y)$ if $x ~ r_i$ and $y~r_j$ for $i
%\neq j$. This metric is semi-algebraic by Tarski-Seidenberg principle
% and it is easy to see that it defines the classical
%topology. \proofend \\
\emph{Proof.} Let $d$ be the metric from corollary \ref{Met}. If a
set $A$ is closed in the classical topology, then the distance
$d(x,A):=\inf \limits _{y\in A}d(x,y)$ is strictly positive for
all points $x$ outside $A$. Now define $F_i:M \rightarrow \R$ by
$F_i(x)=d(x,M\setminus U_i)$. It is semi-algebraic by the
Tarski-Seidenberg principle. Define $G= (\sum \limits _{i=1}^n
F_i)/2n$ and $V_i=\{x\in M | \quad F_i(x)>G(x)\}$. It is easy to
see that $V_i$ is a proper refinement of $U_i$. \proofend \\
In order to formulate the central proposition of this section, we
need to define one technical notion.
\begin{definition}
A collection of continuous semi-algebraic functions $\{F_i\}$ is
called \textbf{basic collection} if every one of them is basic,
and in every point of $M$ one of them is larger than 1.
\end{definition}
%
%\begin{definition}
%Let $M=\cup \limits _{i=1}^n U_i$ be a finite open
%(semi-algebraic) cover of a Nash manifold $M$. A collection of
%continuous semi-algebraic functions $\{F_j\}$ is called
%\textbf{refinement} of the cover $\{U_i\}$ if for any $i$ there
%exists $j$ such that $M_{F_j} \subset U_i$.
%\end{definition}
%
%\begin{proposition} \label{FinCov}
%Let $M=\bigcup \limits _{i=1}^n U_i$ be a finite open
%(semi-algebraic) cover of a Nash manifold $M$. Then there exist
%finite open (semi-algebraic) covers $M=\bigcup \limits _{j=1}^m
%V_j$, $M=\bigcup \limits  _{j=1}^m V'_j$ and $M=\bigcup \limits
%_{j=1}^m V''_j$ and continuous semi-algebraic functions $f_j:M
%\rightarrow \R^{\geq 0}$ such that $V_j$ are affine, the cover
%$V_j$ is a refinement of $U_i$, the cover $V'_j$ is a proper
%refinement of the cover $V_j$, $V'_j=\{x\in V_j|f_j(x) \neq 0\}$,
%$f_j|_{V_j'}$ is Nash and $V''_j=\{x \in V_j|f_j(x)>1\}$.
%\end{proposition}
%
\begin{proposition} \label{FinCov}
Let $M=\bigcup \limits _{i=1}^n U_i$ be a finite open
(semi-algebraic) cover of a Nash manifold $M$. Then there exists a
basic collection of continuous semi-algebraic functions $F_j$ on
$M$ such that the cover $M_{F_j}$ is a refinement of $\{U_i\}$.
\end{proposition}

\emph{Proof.} Cover $U_i$ by affine open subsets $V_{ij}$. By
proposition \ref{NarCov} there exist $V'_{ij} \subset
\overline{V'_{ij}} \subset V_{ij}$ which also cover $M$. By
proposition \ref{FinCor} which follows from the finiteness
theorem, there exist functions $G_{ijk}:V_{ij} \rightarrow \R$ and
open sets $V'_{ijk}$ such that $V'_{ijk}= \{x \in V_{ij} |
G_{ijk}(x) \neq 0\}$, $G_{ijk}|_{V'_{ijk}}$ is positive and Nash
and $\bigcup \limits _k V'_{ijk}=V'_{ij}$. In order to have a
unified system of indexes we denote $V_{ijk}:=V_{ij}$. It gives a
finite cover of $M$ which is a refinement of $U_i$, we re-index it
to one index cover $V_l$. By the same re-indexation we define
$G_l$ and $V'_l$. Extend $G_l$ by zero to a function
$\widetilde{G_l}$ on $M$. It is continuous. Denote $G= (\sum
\widetilde{G_l})/(2n)$ where $n$ is the number of values of the
index $l$. Consider $G|_{V_l}$. This is a strictly positive
continuous semi-algebraic function on an affine Nash manifold,
hence by proposition \ref{MajNash} it can be bounded from below by
a strictly positive Nash function $g'_l$. Denote
$H_l:={G_l}/g'_l$. Extending $H_l$ by zero outside $V_l$ we obtain
a collection of continuous semi-algebraic functions $F_l$ to $M$.
It is easy to see that $\{F_l\}$ is a basic collection and
$M_{F_l}$ is a refinement of $U_i$. \proofend
\begin{remark}
We do not know whether the sets $M_F$, where $F$ is a basic
function, form a basis of the restricted topology for a non-affine
Nash manifold $M$. May be this statement is difficult to verify
for the following reason: the notion of basis is not an
appropriate notion for restricted topology. The property we have
just proven is slightly weaker, but enough for our purposes and
probably can be formulated for any Grothendieck site, unlike the
usual notion of basis.
\end{remark}

\subsection{Cosheaves over restricted topological spaces} \label{cosheaf}

Since Schwartz functions cannot be restricted to open subsets, but
can be extended by 0 from open subsets, we need the notions of a
pre-cosheaf and a cosheaf. There are such notions for any
Grothendieck site (namely, a pre-sheaf and a sheaf with values in
the opposite category). We will now repeat their definitions for
restricted topology.
\begin{definition}
{A \textbf{pre-cosheaf $F$} on a restricted topological space $M$
is a covariant functor from the category ($Top(M)$ which has open
sets as its objects and inclusions as morphisms) to the category
of abelian groups, vector spaces, etc.

In other words, it is an assignment $U \mapsto F(U)$ for every
open $U$ of an abelian group, vector space, etc., and for every
inclusion of open sets $V \subset U$ an extension morphism
$ext_{V,U}: F(V) \rightarrow F(U)$ which satisfy: $ext_{U,U} = Id$
and for $W \subset V \subset U$, $ext_{V,U} \circ ext_{W,V} =
ext_{W,U}$.}
\end{definition}
\begin{definition}
{A \textbf{cosheaf $F$} on a restricted topological space $M$ is a
pre-cosheaf on $M$ fulfilling the conditions dual to the usual
sheaf conditions, and with only finite open covers allowed. This
means: for any open set $U$ and any finite cover $U_i$ of $M$ by
open subsets, the sequence $$ \bigoplus_{i=1}^{n-1}
\bigoplus_{j=i+1}^{n} F(U_i \cap U_j) \rightarrow
\bigoplus_{i=1}^n F(U_i) \rightarrow F(U) \rightarrow 0$$ is
exact. \\
Here, the first map is defined by
$$ \bigoplus_{i=1}^{n-1} \bigoplus_{j=i+1}^{n} \xi_{ij} \mapsto \sum_{i=1}^{n-1} \sum_{j=i+1}^{n}
 ext_{U_i \cap U_j,U_i}(\xi_{ij}) -
 ext_{U_i \cap U_j,U_j}(\xi_{ij}) $$
and the second one by $$\bigoplus_{i=1}^n \xi_i \mapsto \sum
\limits  _{i=1}^n ext_{U_i,U}(\xi_i).$$  }
\end{definition}

\end{document}